\documentclass[reqno,12pt] {article}
\usepackage{amsmath, amsfonts, amssymb,amsthm,hyperref,floatrow}
\usepackage{graphicx, color,dsfont,epsfig,caption,wrapfig,subfig}
\usepackage{xcolor}
\usepackage{tikz}

\hypersetup{
    linktoc=page,
    linkcolor=red,          
    citecolor=blue,        
    filecolor=blue,      
    urlcolor=cyan,
    colorlinks=true           
}

\newtheorem{theorem}{Theorem}
\newtheorem{remark}[theorem]{Remark}
\newtheorem{lemma}[theorem]{Lemma}
\newtheorem{proposition}[theorem]{Proposition}

\newtheorem{corollary}[theorem]{Corollary}
\newtheorem{definition}[theorem]{Definition}

\def \P{ \mathbb P  }
\def \E{ \mathbb E  }

\definecolor{remi}{rgb}{1,0,0}
\usepackage{titlesec}
    \titleformat{\section}[hang]
        {\color{remi}{}\bfseries\filcenter\large}
        {\thesection.}
        {0.4em}
        {}[]

\DeclareMathSymbol{\leqslant}{\mathalpha}{AMSa}{"36} 
\DeclareMathSymbol{\geqslant}{\mathalpha}{AMSa}{"3E} 
\DeclareMathSymbol{\eset}{\mathalpha}{AMSb}{"3F}     
\renewcommand{\leq}{\;\leqslant\;}                   
\renewcommand{\geq}{\;\geqslant\;}                   

\newcommand{\C}{\mathbb{C}}
\newcommand{\R}{\mathbb{R}}

\newcommand{\N}{\mathbb{N}}

\def \Cv{\mathrm{Cov}}
\newcommand{\ind}{\mathds{1}}

\begin{document}

\title{Gaussian multiplicative chaos and KPZ duality}

\author{Julien Barral \footnote{Universit{\'e} Paris 13, Institut Galil\'ee, LAGA, UMR CNRS 7539, 99 rue Jean-Baptiste Cl\'ement 93430 Villetaneuse, France.}, Xiong Jin\footnote{University of St Andrews, Mathematics Institute, North Haugh, St Andrews,
KY16 9SS, Scotland.}, R\'emi Rhodes \footnote{Universit{\'e} Paris-Dauphine, Ceremade, UMR 7534, Place du marechal de Lattre de Tassigny, 75775 Paris Cedex 16, France.} , Vincent Vargas  \footnote{Universit{\'e} Paris-Dauphine, Ceremade, UMR 7534, Place du marechal de Lattre de Tassigny, 75775 Paris Cedex 16, France.} \footnote{The last two authors are partially supported by the CHAMU project (ANR-11-JCJC).}
}
\maketitle

\begin{abstract} This paper is concerned with the construction of atomic Gaussian multiplicative chaos and the KPZ formula in Liouville quantum gravity. On the first hand, we construct purely atomic random measures corresponding to values of the parameter $\gamma^2$ beyond the transition phase (i.e. $\gamma^2>2d$) and check the duality relation with sub-critical Gaussian multiplicative chaos.  On the other hand, we give a simplified  proof of the classical KPZ formula as well as the dual  KPZ formula for atomic Gaussian multiplicative chaos. In particular, this framework allows to construct singular Liouville measures and  to understand the duality relation in Liouville quantum gravity.
\end{abstract}
\tableofcontents

\normalsize

\section{Introduction}
Log-normal multiplicative martingales were introduced by Mandelbrot \cite{mandelbrot1} in order to  build random measures describing energy dissipation and contribute explaining intermittency effects in  Kolmogorov's theory of fully developed turbulence (see \cite{cf:Castaing,cf:Sch,cf:Sto,cf:Cas,cf:Fr} and references therein). However, his model was difficult to be mathematically founded in complete rigor and  this is why he proposed  in \cite{mandelbrot} the simpler model of random multiplicative cascades whose detailed study started with Kahane's and Peyri\`ere's notes~\cite{kahane74,pey74}, improved and gathered in their joint paper~\cite{KP}.

From that moment on, multiplicative cascades have been widely used as reference models in many applications. However, they possess many drawbacks related to their discrete scale invariance, mainly they involve a particular scale ratio and they do not possess stationary fluctuations (this comes from the fact that they are constructed on a $p$-adic tree structure).
In the eighties, Kahane \cite{Kah} came back  to Mandelbrot's initial model and developed 
a continuous parameter theory of suitable multifractal random measures, called Gaussian multiplicative chaos. His efforts were followed by several authors \cite{Bar,cf:Sch,bacry,cf:RoVa,rhovar,Fan,allez,sohier} coming up with various generalizations at different scales. This family of random fields  has found many applications in various fields of science, especially in turbulence and  in mathematical finance. Briefly, a Gaussian multiplicative chaos can be formally understood as a random measure on the Borelian subsets of (a domain of) $\R^d$ of the form
\begin{equation}\label{eq:introchaos1}
M(A)=\int_Ae^{X_x-\frac{1}{2}\E[X_x^2]}\,dx,
\end{equation}
where $(X_x)_x$ is a centered Gaussian field and $dx$ is the Lebesgue measure. In his seminal work \cite{Kah}, Kahane showed that the only short scale behaviour of the correlations of the field $X$ that produces interesting measures is logarithmic.  In the case of logarithmic correlations, it is standard to rewrite equation \eqref{eq:introchaos1} as
\begin{equation}\label{eq:introchaos2}
M_\gamma(A)=\int_Ae^{\gamma X_x-\frac{\gamma^2}{2}\E[X_x^2]}\,dx,
\end{equation}
where $\gamma>0$ is a parameter, sometimes called intermittency parameter or coupling constant depending on the context, and $X$ is a centered Gaussian field with covariance of the form:
\begin{equation}\label{eq:correl}
\E[X_xX_y]=\ln_+\frac{1}{|x-y|}+g(x,y)
\end{equation}
for some continuous bounded function $g$ (and $\ln_+x=\max(\ln x,0)$). Interestingly and in this context, these measures present a phase transition: they are non trivial if and only if $\gamma^2<2d$ (see \cite{Kah}). Let us also mention that in \cite{Kah}, Kahane made a thorough study of these measures (existence of moments, non degeneracy,  etc...).

Recently and in a major conceptual step, the authors in \cite{cf:DuSh} have drawn attention on the fact that 2d-Gaussian multiplicative chaos should be considered to give a rigorous meaning to the Liouville measure in Liouville Quantum Gravity (see \cite{cf:KPZ,cf:Da,cf:DuSh,bourbaki} among many others). Based on \cite{cf:KPZ,cf:Da}, they argue that the field $X$ in \eqref{eq:introchaos1} has to be a Gaussian Free Field (GFF for short) to obtain the Liouville measure (see \cite{cf:DuSh,Kah,cf:RoVa} for the construction of the measures). In this context, the KPZ formula has been proved rigorously \cite{cf:DuSh,cf:RhoVar} (see also \cite{Benj} in the context of multiplicative cascades) below the phase transition arising at $\gamma^2=4$ (i.e. $\gamma^2=2d$ with $d=2$), where the constant $\gamma$ is related to the central charge $c\leq 1$ of the underlying conformal field theory by the relation (see \cite{cf:KPZ})
$$ c=1-\frac{4}{6}\big(\gamma-\frac{4}{\gamma}\big)^2.$$ 
Actually, the above relation produces two different solutions (for $c<1$), a first one in $\gamma\in [0,2[$ and a second one $\bar{\gamma}\in]2,+\infty[$, both of them satisfying the relation $\gamma\bar{\gamma}=4$. The mathematical investigations in \cite{Benj,cf:DuSh,cf:RhoVar} (based on a huge amount of physics papers) deal with the $\gamma<2$ solution, called the {\it standard branch of gravity}. Physicists have also investigated the implications of the existence of the other solution $\bar{\gamma}\in]2,+\infty[$, starting with the works \cite{Kleb1,Kleb2,Kleb3}, in what they called the {\it non-standard branch of gravity} or the {\it dual branch of gravity}. It is clear that Kahane's theory does not allow to go beyond the $\gamma^2=4$ threshold so that one has to look for other types of measures than those considered in \cite{Kah,cf:DuSh,cf:RhoVar} to explain the {\it dual branch of gravity}, i.e. to construct mathematically singular  Liouville measures beyond the  phase  transition (i.e. for $\gamma^2>4$). 

If one considers Mandelbrot's multiplicative cascades as a powerful toy model to understand the continuous model, or Gaussian multiplicative chaos, then one has to go back to Durrett and Liggett's deep paper \cite{durrett} to have a rather good intuitive picture. The authors solve in full generality an equation, called Mandelbrot's star equation, which is satisfied by the law of the total mass of the  limit of a non degenerate Mandelbrot multiplicative cascade, and describes the scale invariance of such a measure. However this equation captures the scale invariance of a larger class of random measures, among which limit of Mandelbrot cascades are characterized as being those with total mass in $L^1$. Sticking to our notations, when $\bar{\gamma}>2$, one can deduce from \cite{durrett}  that  the random measure associated with the star equation can be nothing but the derivative (in the distribution sense) of  a stable subordinator subordinated by a sub-critical multiplicative cascade ($\gamma<2$) or a stable subordinator subordinated by the critical multiplicative cascade ($\gamma=2$). We do not detail here what the critical case is but we refer the reader to \cite{DRSV,DRSV2} for recent results in the topic of Gaussian multiplicative chaos and an extensive list of references about the origins of this concept in various related models. A continuous analog of the Mandelbrot's star equation, called star scale invariance, has been introduced in \cite{allez} so that one is left with the intuition that the same picture can be drawn for Gaussian multiplicative chaos. Among the two options we are left with (subordinating by a sub-critical Gaussian multiplicative chaos or a critical Gaussian multiplicative chaos), it turns out that subordinating by a critical Gaussian multiplicative chaos cannot be a suitable model for duality. Indeed, such measures do not satisfy the expected KPZ relations (see below). So, duality must necessarily correspond to subordinating by a sub-critical Gaussian multiplicative chaos. We also stress that the notion of star scale invariance is nothing but a rigorous formulation of the scaling heuristics developed in \cite{PRL} to quantify the measure of a Euclidean ball of size $\epsilon$ (see in particular the section \textit{Liouville quantum duality}).  Let us further stress that  the concept of subordinating sub-critical Gaussian multiplicative chaos by L\'evy processes has been studied for the new and fundamental examples of multifractal formalism illustrations it provides \cite{BS}. This general procedure consisting in subordinating  self-similar processes by nice multifractal measures to create new models of multifractal processes has been pointed out by Mandelbrot, in particular in \cite{Mfin}. After posting several versions of this paper online, we learnt the existence of \cite{Dup:houches} where the author suggests to use the subordination procedure  to model the dual branch of gravity, confirming in a way the above heuristics. As pointed out in \cite{Dup:houches}, the subordination procedure in dimension higher than $1$ requires  a small additional mathematical machinery. Indeed, while in dimension $1$ it is enough to subordinate a standard stable subordinator by an independent Gaussian multiplicative chaos \cite{BS,durrett}, in higher dimensions the canonical construction is to define a stable independently scattered random measure \cite{Ros} conditionally  on a random intensity, which is given by a Gaussian multiplicative chaos. We give here a simple argument to prove that this can be achieved in a measurable way, thus justifying the fact that the measure proposed in \cite{Dup:houches} is mathematically founded (see Remark \ref{rem.otherlaw}). 

Another natural construction of such measures is to let  a multiplicative chaos act on independent atomic measures: this is the first main issue of this paper. Though producing the same measure in law as the subordination procedure, it is conceptually more powerful and flexible. To motivate this comment, we first draw up the framework a bit more precisely and we will come back to this point thereafter. Fix a simply connected domain $D\subset \C$. For $\gamma <2$, the {\it standard Liouville measure} can formally be written as 
\begin{equation}\label{measintro}
M_\gamma(A)=\int_Ae^{\gamma X_x-\frac{\gamma^2}{2}\E[X_x^2]}\,dx
\end{equation}
 where $X$ is the Gaussian Free Field (GFF) over the domain $D$. In fact and to be exhaustive, $M_\gamma(A)=\int_Ae^{\gamma X_x-\frac{\gamma^2}{2}\E[X_x^2]}\,h(x)dx$ where $h(x)$ is a deterministic function defined in terms of  the conformal radius, but this additional term $h$ does not play a role in the nature of the problems/results we want to address here. For a given compact set $K\subset D$, it has been proved that the Hausdorff dimension of $K$ computed with the Euclidian metric, call it ${\rm dim}_{Leb}(K)$, is related to the Hausdorff dimension of $K$ computed with the measure $M_\gamma$, call it ${\rm dim}_{ \gamma}(K)$. The connection is the so-called KPZ formula
$${\rm dim}_{Leb}(K)=(1+\frac{\gamma^2}{4}){\rm dim}_{\gamma}(K)-\frac{\gamma^2}{4} {\rm dim}_{\gamma}(K)^2.$$
The measure based approach of duality in Liouville quantum gravity is about constructing  purely atomic random measures $M_{\bar{\gamma}}$, for parameter values $\bar{\gamma}>2$ (i.e. beyond the phase  transition), that satisfy the KPZ relation
 $${\rm dim}_{Leb}(K)=(1+\frac{\bar{\gamma}^2}{4}){\rm dim}_{\bar{\gamma}}(K)-\frac{\bar{\gamma}^2}{4} {\rm dim}_{\bar{\gamma}}(K)^2.$$
Then, by considering the dual value $\gamma=\frac{4}{\bar{\gamma}}$ of the parameter $\bar{\gamma}$, one wants to rigorously establish the duality relation
$${\rm dim}_{\bar{\gamma}}(K)=\frac{\gamma^2}{4}{\rm dim}_{ \gamma}(K) .$$
We point out that physicists can recover the (more classical) relation between the scaling exponents by setting $\triangle_\gamma=1-{\rm dim}_{ \gamma}(K) $ and $\triangle_{\bar{\gamma}}=1-{\rm dim}_{ \bar{\gamma}}(K) $. Proving rigorously the dual KPZ formula is the second main point of our paper. Let us also mention that scaling heuristics have been developed in \cite{Dup:houches,PRL}. It is not straithforward to derive a rigorous KPZ statement from \cite{Dup:houches,PRL} as it involves working under the Peyri\`ere measure which in the dual context is infinite (the authors of \cite{Dup:houches,PRL} call this measure the rooted measure).
 
Consider $\alpha=\frac{\gamma^2}{4}\in]0,1[$ and define a random measure $M_{\bar{\gamma}}$, the law of which is given, conditionally  to $M_\gamma$, by a independently scattered $\alpha$-stable random measure with control measure $M_\gamma$. More precisely and conditionally  to $M_\gamma$, we consider a Poisson random measure $n_{\gamma,\alpha}$ distributed on $D\times\R_+^*$ with intensity
$$M_\gamma(dx)\frac{dz}{z^{1+\alpha}}.$$
Then the law of the dual measure is given by
\begin{equation}\label{mbar:introsub}
\forall A\subset D \text{ Borelian},\quad M_{\bar{\gamma}}(A)=\int_A\int_0^{+\infty}z\,n_{\gamma,\alpha}(dx,dz).
\end{equation}
This is the subordination procedure proposed in \cite{Dup:houches}. Nevertheless, we feel that this construction suffers a conceptual drawback: there is no way to obtain this measure as the almost sure limit of a properly regularized sequence, meaning of a sequence where the GFF distribution has been suitably regularized via a cutoff procedure.  This point is important. Indeed, if one has a look at the original physics motivations \cite{Kleb1,Kleb2,Kleb3}, dual gravity should be produced by ``pinching" suitably the random field under consideration: in \cite{Kleb1,Kleb2,Kleb3}, the terminology ``finely tuned" is used in the context of random matrix theory.  In terms of planar maps, this should correspond to creating enhanced bottlenecks in the planar map \cite{Kleb1}. Therefore, it is natural to think that this pinching will have the effect of creating a bi-randomness: production of atoms on the one hand and creation of a multiplicative chaos term on the other hand. These two terms are likely to mix together to produce a random measure, the law of which is given by \eqref{mbar:introsub}. For instance, we roughly explain a scenario that may produces such measures. Consider a random measure on a bounded open subset of $\R^2$ with density given by the exponential of random field. Assume that this field is dominated by two independent scales: a first one with strong local interactions like a (regularized) GFF and a second independent scale with a stronger independence structure which, for instance, macroscopically behaves like a Random Energy Model (REM for short). If the temperature of the REM scale is low enough, the REM-like scale  component will macroscopically behave like an independent stable random measure. Therefore, after suitable renormalization, the random measure is well approximated by the exponential of a GFF integrated by an independent stable random measure. We will call these measures {\it atomic Gaussian multiplicative chaos}. Observe that the location of the atoms is here independent of the GFF, unlike the subordination procedure where position of atoms is strongly coupled to the realization of the underlying Gaussian multiplicative chaos. It is therefore interesting to study the convergence of such a structure, if possible in a strong sense like almost sure convergence  or convergence in probability. We further stress that this point has never been investigated either in the whole literature of Gaussian multiplicative chaos or in the literature concerning closely related models: branching random walks, multiplicative cascades or branching Brownian motions.
 
The rough picture we have just drawn can be mathematically formulated as follows. Consider a couple of exponents $(\gamma,\bar{\gamma})$ such that $\gamma^2<4$ and $\gamma\bar{\gamma}=4$ (this relation between the exponents may be explained via a volume preserving condition, see Section \ref{sec:KPZ}). Introduce an independently scattered random measure $n_\alpha$ characterized by its Laplace transform ($|A|$ stands for the Lebesgue measure of $A$)
$$ \forall A \subset \R^2 \text{ Borelian}, \quad \E[e^{-u n_\alpha(A)}]=e^{-u^\alpha|A|}$$
where $\alpha=\frac{\gamma^2}{4}$. Consider a suitable cutoff approximation family $(X^n)_n$ of the GFF distribution (see section \ref{back} for a discussion about possible cutoff procedures). Prove that the family of approximated measures
\begin{equation}\label{dualintron}
M^n_{\bar{\gamma}}(A)=\int_Ae^{\bar{\gamma} X^n_x-2\E[(X^n_x)^2]}\,n_\alpha(dx),
\end{equation}
converges almost surely (possibly along a subsequence) towards a measure that we will formally write as:
\begin{equation}\label{dualintro}
M_{\bar{\gamma}}(A)=\int_Ae^{\bar{\gamma} X_x-2\E[X_x^2]}\,n_\alpha(dx).
\end{equation}
We point out that the above atomic Gaussian multiplicative chaos is not a Gaussian multiplicative chaos in the usual sense. Indeed the lognormal weight is not normalized to have expectation $1$ and the above renormalization is therefore not standard. Actually, the expectation explodes giving rise to a strong competition between the atoms produced by the random measure $n_\alpha$ and the ability of the lognormal weight to kill these atoms. 

We will prove that the measure $(M^n_{\bar{\gamma}})_n$  converges almost surely along deterministic subsequences, giving sense to a new and exciting theory of (non-standard) multiplicative chaos with respect to atomic measures. It may be worth mentioning here that this construction produces, as expected from \cite{durrett}, star scale invariant random measures for suitable choices of the distribution $X$. Beyond the applications in Liouville Quantum Gravity, we have the feeling that this approach offers new perspectives in the theory of Gaussian multiplicative chaos  that we develop in Section \ref{pers}. Therefore we will present our construction in a more general context than the only purpose of giving sense to the {\it dual branch of gravity}: we will not restrict ourselves to the $2$-dimensional case and we will not consider the only GFF but, more generally, log-correlated Gaussian distributions.
 
 \section{Background}\label{back}
 In this section, we will briefly explain Kahane's theory of multiplicative chaos in $\R^d$. In fact, Kahane's theory is valid in any open domain $D \subset \R^d$ with no substantial change. At the end of the section, we will also roughly recall the connection with measures formally given by the exponential of the GFF. 
 
\subsection{Sigma positive kernels}\label{back:kah}
We consider a covariance kernel $K$ of $\sigma$-positive type (\cite{Kah}), namely that  $K$ can be rewritten as a sum
\begin{equation}\label{decompo}
\forall x,y\in\R^d,\quad K(x,y)=\sum_{n\geq 1}q_n(x,y)
\end{equation} where $(q_n)_n$ is a sequence of continuous positive kernels of positive type. We further assume that 
\begin{equation}\label{formgam}
\forall x\in\R^d,\quad K(x,y)= \ln_+\frac{T}{|x-y|}+g(x,y)
\end{equation}
 where $g$ is a bounded continuous function over $\R^d\times \R^d$ (and $\ln_+(x)=\max(0,\ln(x))$). We can consider a sequence of independent centered Gaussian processes $(Y^n)_{n\geq 1} $ where, for each $n\geq 1$,  $(Y^n_x)_{x\in\R^d}$ is a centered continuous Gaussian field with covariance function given by
$$\forall x,y\in \R^d,\quad \Cv(Y^n_x,Y^n_y)=q_n(x,y).$$
Finally, for $n\geq 1$, we define:
$$X^n_x=\sum_{p=1}^nY^p_x.$$
It is a centered continuous Gaussian process with covariance function:
\begin{equation}
\forall x,y\in \R^d,\quad k_n(x,y)\stackrel{def}{=}\Cv(X^n_x,X^n_y)=\sum_{k=1}^nq_k(x,y).
\end{equation}
The reader may find several important examples of sigma-positive kernels in Appendix \ref{rem:pas}. We stress that the Gaussian distribution $X$ is not required to be stationary.
 
\subsection{Gaussian multiplicative chaos}
For each $n\geq 1$, we can define a Radon measure $M_n$ on the Borelian subsets of $\R^d$ by
\begin{equation}\label{defM}
M_n(A)=\int_Ae^{\gamma X^n_x-\frac{\gamma^2}{2}\E[(X^n_x)^2]}\,dx.
\end{equation}
For each Borelian set $A$, the sequence $(M^n(A))_n$ is a positive martingale. Thus it converges almost surely towards a random variable denoted by $M(A)$. One can deduce that the sequence of measures $(M_n)_n$ weakly converges towards a random Radon measure $M$, commonly denoted by 
\begin{equation}\label{GMC}
M(A)=\int_Ae^{\gamma X_x-\frac{\gamma^2}{2}\E[X_x^2]}\,dx 
\end{equation}
and called Gaussian multiplicative chaos associated to the kernel $\gamma^2 K$. Roughly speaking, (\ref{GMC}) can be understood as a measure admitting as density the exponential of a Gaussian distribution $X$ with covariance kernel $\gamma^2 K$. Of course, this is purely formal because the exponential of a random distribution cannot be directly defined. Kahane proved that the martingale $(M_n(A))_n$, for some Borelian set $A$ with non-null finite Lebesgue measure, is uniformly integrable if and only if $\gamma^2<2d$, in which case it has no atoms. This condition $\gamma^2<2d$ is necessary and sufficient in order for the limiting measure $M$ to be non identically null. Furthermore, he proved that the law of the limiting measure $M$ does not depend on the decomposition \eqref{decompo}
of $K$ into a sum of positive continuous kernels, hence founding an intrisic theory. 

When  $\gamma^2<2d$, the power-law spectrum $\xi$, defined through the relation
$$\E[M(B(0,\lambda)^q]\simeq C_q \lambda^{\xi(q)},\quad \lambda\to 0$$ for all $q\geq 0$ such that the expectation makes sense (i.e. for $0\leq q <\frac{2d}{\gamma^2}$, see \cite{Kah}),
 is given by
$$\xi(q)=(d+\frac{\gamma^2}{2})q-\frac{\gamma^2}{2}q^2.$$

\subsection{Generalized Gaussian multiplicative chaos}

In a series of papers \cite{cf:RoVa1,cf:RoVa}, the authors developped a generalized theory of Gaussian multiplicative chaos in all dimensions and for all translation invariant logarithmic kernels of positive type. They prove that one can construct the measure as the limit in law of measures where the underlying Gaussian field has a covariance kernel obtained by a smooth cutoff of the logarithmic one (their approach relies on convolutions of the covariance function rather than on the notion of $\sigma$-positive type (\cite{Kah}). In particular, they show that the law of the limit measure does not depend on the smoothing procedure one applies to the logarithmic kernel. Their techniques are quite general and work with any smooth cutoff of the covariance kernel (or the underlying field).

\subsection{Application to the construction of Liouville measures in dimension $2$}\label{sub:GFF}
Formally, the GFF (or Euclidian bosonic massless free field) in a bounded domain $D\subset \R^2$ is a ``Gaussian Field" $X$ with covariance given by:
\begin{equation*}
\E[X_xX_y]=G(x,y),
\end{equation*}
where $G $ is the Green function of $D$ with zero boundary condition (see for instance \cite{ShGFF} or chapter 2.4 in \cite{cf:Law} for the definition and main properties). Let $B$ be a Brownian motion starting from
 $x\in D$ under the measure $P^x$ and consider the stopping time $T_D=\text{inf} \{t \geq 0, \; B_t\not \in D \}$.
  If we denote $p_D(t,x,y)=P^x(B_{t} \in dy, \; T_D > t)$, we have:
\begin{equation*}
G (x,y)= \pi \int_{0}^{\infty}p_D(t,x,y)dt.
 \end{equation*}
Note that, for each $t>0$, $p_D(t,x,y)$ is a continuous positive and positive definite kernel on $D$. Therefore, following Kahane's theory, we can define the Gaussian multiplicative chaos $M$ associated to the
  kernel $\gamma^2G$. Since the Green function takes on the form \eqref{formgam}, this measure is not trivial provided that $\gamma^2<4$.

After the two theories mentioned above (Gaussian multiplicative chaos and Generalized Gaussian multiplicative chaos),  the authors of  \cite{cf:DuSh} later suggested two other constructions of the Liouville measure. The first construction, which falls under the scope of Kahane's theory, consists in expanding the Green function of the Laplacian along an orthonormal basis of the Sobolev space $H^1_0(D)$ (the space of functions whose gradient is in $L^2(D)$). All the results presented in this paper applies to this construction.  The second construction consists in averaging the GFF distribution along circles of size $\epsilon$ to get an $\epsilon$-regularized approximation of the GFF. By the machinery developped in the generalized Gaussian multiplicative chaos theory \cite{cf:RoVa1,cf:RoVa}, it is straightforward to see that the measure converges in law to Gaussian multiplicative chaos with the Green function as covariance kernel (in fact this is true for ball averages or any smooth cutoff procedure). The contribution of the authors in \cite{cf:DuSh} is to prove almost sure convergence along powers of 2 in the case of circle average approximations and to prove almost sure equality of the measures obtained via  circle average/$H^1_0(D)$-expansion approximations of the GFF.

\section{Atomic Gaussian multiplicative chaos}\label{sec:chaos}

\begin{remark} We stick to the notations of the previous section. We nevertheless assume that the considered Gaussian fields are stationary. Though it may appear as a restriction, the proofs in the general case work exactly the same. Actually, being stationary or not is just hidden in the ``small noise $g$" appearing in \eqref{formgam}. \end{remark}
 
Now we begin to construct what will be called atomic Gaussian multiplicative chaos. To this purpose, we consider a Gaussian multiplicative chaos, denoted by $M$,  for some $\gamma^2<2d$. Choose $\alpha\in ]0,1[$ and  consider a  Poisson random measure $N_\alpha$ distributed on $\R^d\times \R^*_+$ with intensity $dx\frac{dz}{z^{1+\alpha}} $ and independent of the sequence $(Y^n_x)_{x\in\R^d}$. We introduce the random measure $$n_\alpha(dx)=\int_0^{+\infty}z\,N_\alpha(dx,dz) ,$$ which can be thought of as an independently scattered stable random measure.
Then, for $\gamma^2<2d$, we define the sequence of random measures
\begin{align}
\forall A\in \mathcal{B}(\R^d),\quad \overline{M}_n(A)=&\int_Ae^{\frac{\gamma}{\alpha} X^n_x-\frac{\gamma^2}{2\alpha}\E[(X^n_x)^2]}n_\alpha(dx).\label{defMbar2}
\end{align}

\begin{theorem}\label{as}{(Convergence in probability)}\\
(1) For each bounded Borelian set $A$, the sequence $(\overline{M}_n(A))_n$ converges in probability towards a non trivial random variable. \\
(2) For each subsequence, we can extract a (deterministic) subsequence such that, almost surely, the sequence of random measures $(\overline{M}_n(dx))_n$ weakly converges  towards a random Radon measure~$\overline{M}$.\\
(3) the law of $\overline{M}$ is characterized by the following relation:
 \begin{equation}\label{relfundn}
\E[e^{-u_1 \overline{M}(  A_1)-\dots-u_p\overline{M}(  A_p)}]=\E\big[ e^{-\frac{\Gamma(1-\alpha)}{\alpha}\big(u_1^\alpha M(A_1)+\dots+u_p^\alpha M(A_p)\big)} \big]
\end{equation} valid for all $u_1,\dots,u_p\in\R_+$ and all disjoint Borelian subsets $A_1,\dots,A_p\subset \R^d$. \\
(4) the limiting measure  $\overline{M}$ is non trivial for $\gamma^2<2d$ and all $\alpha\in]0,1[$. 
\end{theorem}

Observe that formula \eqref{relfundn} is the standard Levy-Khintchine relation for independently scattered random measures (see \cite{Ros}).  

\begin{definition}{\bf Atomic Gaussian multiplicative chaos}
The limit random measure $\overline{M}$ defined  in Theorem \ref{as} will be called atomic Gaussian multiplicative chaos. It will be formally written as  (with $\overline{\gamma} =\frac{\gamma}{\alpha}$)
\begin{equation}\label{formal}
\overline{M}(\cdot)=\int_\cdot e^{\overline{\gamma}  X_x-\frac{\alpha\overline{\gamma}^2}{2}\E[X_x^2]}n_\alpha(dx)
\end{equation}
where $X$ is a stationary Gaussian distribution with covariance kernel $K$. 
\end{definition}
The above expression \eqref{formal}   justifies the fact that the measure $\overline{M}$ can be seen as a non standard Gaussian multiplicative chaos since the weight has not expectation $1$. Furthermore, it can be defined for values of $\overline{\gamma}^2$ beyond the critical value $\overline{\gamma}^2=2d$. Notice that the renormalization (i.e. $\frac{\alpha\overline{\gamma}^2}{2}\E[X_x^2]$) differs from the standard Gaussian multiplicative chaos.

\begin{proposition}\label{def}
For $\gamma^2<2d$ and $\alpha\in]0,1[$, the law of the random measure $\overline{M}$ does not depend on the decomposition of $K$ into a sum of positive continuous kernels of positive type. Furthermore, $\overline{M}$ is almost surely a purely atomic measure.
\end{proposition}

\begin{remark}\label{rem.otherlaw}
 There is another way of seeing the law of the measure $\overline{M}$. We want to introduce a positive Radon random measure $N_M$ distributed on $\R^d\times\R^*_+$, whose law conditionally to $M$ is that of a Poisson random measure with intensity
$$\frac{M(dx)\,dz}{z^{1+\alpha}}.$$ Then we want to consider the family of purely atomic positive random measures
\begin{align} 
\forall A\in \mathcal{B}(\R^d),\quad \widetilde{M}(A)=&\int_A\int_{\R_+}z\,N_M(dx,dz).\label{defMbar}
\end{align}
If one can give a rigorous meaning to the above contruction, then the law of the random measure $\widetilde{M}$ is the same as that of Theorem \ref{as}. Nevertheless, the reader may observe that it is not obvious to give a ``measurable" construction of the above measure $\widetilde{M}$. Of course, Theorem \ref{as} provides a first rigorous way of defining a random measure whose law is given by \eqref{relfundn}. We give in Appendix \ref{otherlaw} an alternative construction.
\end{remark}

\subsection{Power-law spectrum and moments of the atomic chaos}
In this subsection, we assume that $\gamma^2<2d$. Let us define  
$$\forall q\in \R,\quad \overline{\xi}(q)=\big(\frac{d}{\alpha}+\frac{\gamma^2}{2\alpha}\big)q-\frac{\gamma^2}{2\alpha^2}q^2.$$ We will show below that this function coincides with the power law spectrum of the measure $\overline{M}$. In particular, we see that $\overline{\xi}(q)=\xi(\frac{q}{\alpha})$. Furthermore, from Theorem \ref{as}, the Laplace exponents of the measure $\overline{M}$ can be expressed in terms of the Laplace exponents of $M$:
\begin{equation}\label{exp}
 \E[e^{-u \overline{M}(A)}]=\E[e^{-\frac{\Gamma(1-\alpha)}{\alpha}u^\alpha M (A)}].
\end{equation}
This relation allows us to deduce the main properties of the moments of the measure $\overline{M}$:
\begin{proposition}\label{prop:mom}
For all Borelian set $A$ with finite (non null) Lebesgue measure, the random variable $\overline{M}(A)$ possesses a moment of order $\beta\geq 0$ if and only if $\beta<\alpha$. 

Furthermore, we can make explicit the connection between   the moments of $M$ and $\overline{M}$: for all $0\leq \beta <\alpha$,
\begin{equation}\label{relmom}
\E[(\overline{M} (A))^\beta]=   \frac{  \Gamma(1-\beta/\alpha)\Gamma(1-\alpha)^{\beta/\alpha}}{ \Gamma(1- \beta)\alpha^{\beta/\alpha}}  \E[(M(A))^\frac{\beta}{\alpha}]
 \end{equation}
\end{proposition}

\begin{theorem}{\bf (Perfect scaling).}\label{perfect}
If the kernel $K$ is given by 
$$K(x)=\ln_+\frac{T}{|x|}+g(x)$$ where $g$ is a continuous bounded function that is constant in a neighborhood of $0$ then, for some $R>0$:
\begin{equation}\label{scaling}
\forall 0<\lambda<1,\quad (\overline{M}(\lambda A))_{A\subset B(0,R)}\stackrel{law}{=}\lambda^{d/\alpha}e^{\frac{\Omega_\lambda}{\alpha}}(\overline{M}(  A))_{A\subset B(0,R)}
\end{equation}
where $\Omega_\lambda$ is a Gaussian random variable independent of the measure $(\overline{M}(  A))_{A\subset B(0,R)}$ the law of which is characterized by:
$$\E[e^{q \Omega_\lambda}]=\lambda^{\frac{\gamma^2}{2}q-\frac{\gamma^2}{2}q^2}.$$
In particular, for all $0\leq q <\alpha$:
$$\E[\overline{M}(  B(0,\lambda R) )^q]=\lambda^{\overline{\xi}(q)}\E[\overline{M}( B(0,  R) )^q].$$
\end{theorem}

\begin{corollary}\label{notperf}
Assume that the kernel  $K$ takes on the form \eqref{formgam}. Then, for all $0\leq q <\alpha$:
$$\E[\overline{M}(  B(0,\lambda R ) )^q]\simeq C_{q,R} \lambda^{\overline{\xi}(q)}$$ as $\lambda\to 0$ for some positive constant $C_{q,R}$ only depending on $q,R$.
\end{corollary}

\section{KPZ formula and duality}\label{sec:KPZ}
In this section, we adjust the parameters to stick to Liouville quantum gravity issues. Strictly speaking, one should set the dimension equal to $2$ but we will refrain from doing this since the proofs work the same, whatever the dimension is. We consider a standard $\gamma^2<2d$ and its dual exponent $\overline{\gamma} $ by the relation
\begin{equation}\label{relgam}
\gamma \overline{\gamma} =2d.
\end{equation}
This gives $ \overline{\gamma}=\frac{\gamma}{\alpha}$ with
\begin{equation}\label{relgam}
\alpha=\frac{ \gamma^2}{2d}\in]0,1[.
\end{equation}
With this value of $(\overline{\gamma},\alpha)$, the power law spectrum of the corresponding random measure $\overline{M}$ 
\begin{align}
\forall A\in \mathcal{B}(\R^d),\quad \overline{M}(A)=&\int_Ae^{\frac{\gamma}{\alpha} X_x-\frac{\gamma^2}{2\alpha}\E[(X_x)^2]}n_\alpha(dx)\nonumber\\
=&\int_A e^{\overline{\gamma}  X_x-d\E[(X_x)^2]}n_\alpha(dx).\label{defMbar2}
\end{align}
can be rewritten as:
\begin{equation}\label{pwbar}
\overline{\xi}(q)= \big(d+\frac{\overline{\gamma}^2}{2}\big)q-\frac{\overline{\gamma}^2}{2 }q^2.
\end{equation}
Given $\gamma^2<2d$, we stress that this value of $\alpha$ is the only possible value of $0<\alpha<1$ ensuring the (statistically) volume preserving condition $\overline{\xi}(1)=d$. In the language of conformal field theory, this ensures that the theory has no Weyl anomaly (see \cite{cf:Da}).

The KPZ formula is a relation between the Hausdorff dimensions of a given set $A$ as measured by the Lebesgue measure, $M$ or $\overline{M}$. So we first recall how to define these dimensions. Given a Radon measure $\mu$ on $\R^d$ and $s\in [0,1]$, we define 
$$H^{s,\delta}_\mu(A)= \inf \big\{\sum_k \mu(B_k)^{s} \big\}$$ where the infimum runs over all the covering $(B_k)_k$ of $A$ with closed Euclidean balls with radius $r_k\leq \delta$. Clearly, the mapping $\delta>0\mapsto H^{s,\delta}_\mu(A)$ is decreasing. Hence we can define:

$$H^{s}_{\mu}(A)=\lim_{\delta\to 0}H^{s,\delta}_{\mu}(A).$$

$H^{s}_\mu$ is a metric outer measure on $\R^d$ (see \cite{falc} for the definitions). We point out that the fact that $\mu$ possesses atoms or not does not give rise to any additional difficulty.  Thus $H^{s}_\mu$ is a measure on the $\sigma$-field of $H^{s}_\mu$-measurable sets, which contains all the  Borelian sets.

The $\mu$-Hausdorff dimension of the set $A$ is then defined as the value
\begin{equation} 
{\rm dim}_\mu(A)=\inf\{s\geq 0; \,\,H^s_\mu(A)=0\}.
\end{equation}
Notice that ${\rm dim}_\mu(A)\in [0,1]$.
However, it is not clear, in great generality, that we have the classical property:
\begin{equation}\label{HSD}
{\rm dim}_\mu(A)=\sup\{s\geq 0; \,\,H^s_\mu(A)=+\infty\}.
\end{equation}
This is due to the possible presence of atoms for the measure $\mu$.
However we claim

\begin{proposition}\label{prop:HSD}
For any diffuse measure, \eqref{HSD} holds. If the measure $\mu$ possesses atoms, \eqref{HSD} holds for any compact set $A$ that does not encounter the atoms of $\mu$.
\end{proposition}

\begin{corollary}\label{cor:HSD}
If we take $\mu=\text{Leb}$ then \eqref{HSD} holds. If we take $\mu=M$ then, almost surely, \eqref{HSD} holds for every bounded Borelian set. If we take $\mu=\overline{M}$ and $A$ a compact set with null Lebesgue measure then  \eqref{HSD} holds almost surely.
\end{corollary}
This proposition allows to characterize the Hausdorff dimension as the critical value at which 
the mapping $s\mapsto H^s_\mu(A)$ jumps from $+\infty$ to $0$.

In what follows, given a compact set $K$ of $\R^d$ with null Lebesgue measure, we define its Hausdorff dimensions 
${\rm dim}_{Leb}(K)$, ${\rm dim}_M(K)$, ${\rm dim}_{\overline{M}}(K)$  computed as indicated above with $\mu$ respectively equal to the Lebesgue measure, $M$ and $\overline{M}$.

\begin{theorem}\label{KPZ}{{\bf KPZ duality.}}
Let $K$ be a compact set of $\R^d$ with null Lebesgue measure. Almost surely, we have the relations
$${\rm dim}_{Leb}(K)=\frac{\xi({\rm dim}_M(K))}{d}\quad {\rm dim}_{Leb}(K)=\frac{\overline{\xi}({\rm dim}_{\overline{M}}(K))}{d}$$ where $\xi(q)=(d+\frac{\gamma^2}{2})q-\frac{\gamma^2}{2}q^2$ and   $\overline{\xi}(q)= \big(d+\frac{\overline{\gamma}^2}{2}\big)q-\frac{\overline{\gamma}^2}{2 }q^2$. In particular, we have the duality relation between the scaling exponents
\begin{equation}\label{dual}
{\rm dim}_{\overline{M}}(K)= \frac{\gamma^2}{2d}{\rm dim}_{M}(K). 
\end{equation}
\end{theorem}

\vspace{3mm}

\begin{remark}
Note that, in the classical physics literature (in particular $d=2$), it is more usual to focus on the scaling exponents
$$\triangle_{\gamma}=1-{\rm dim}_M(K),\quad  \triangle_{\overline{\gamma}}=1-{\rm dim}_{\overline{M}}(K),\quad x=1-{\rm dim}_{Leb}(K),$$ instead of ${\rm dim}_M(K)$, ${\rm dim}_{\overline{M}}(K)$, ${\rm dim}_{Leb}(K)$. Then the KPZ relations read
$$x=\frac{\gamma^2}{4}\triangle_\gamma^2+(1-\frac{\gamma^2}{4})\triangle_\gamma \quad \text{ and }\quad x=\frac{\overline{\gamma}^2}{4}\triangle_{\overline{\gamma}}^2+(1-\frac{\overline{\gamma}^2}{4})\triangle_{\overline{\gamma}} .$$ The duality relation then becomes
$$\triangle_{\overline{\gamma}}-1 = \frac{\gamma^2}{2d}(\triangle_{\gamma}-1)=\frac{2d}{\overline{\gamma}^2}(\triangle_{\gamma}-1).$$
 \end{remark}

\begin{remark}
If one looks for random measures satisfying the duality relation \eqref{dual}, it is plain to deduce that such a relation implies that the power law spectrum is necessarily given by \eqref{pwbar}. Such a power law spectrum indicates that the searched random measures cannot be defined by \eqref{GMC} in the sense that the integrating measure ($dx$ in \eqref{GMC}) cannot be the Lebesgue measure. Indeed, otherwise Kahane's theory ensures that such measure is identically null. So one has to look for other integrating measures in \eqref{GMC} than the Lebesgue measure. By noticing that necessarily $\bar{\xi}(q)=\xi(\frac{q}{\alpha})$, one can intuitively recover our construction, namely that the searched  measures should be  Gaussian multiplicative chaos integrated against independently scattered $\alpha$-stable random measures, as stated in Theorem \ref{as}.
 \end{remark}

\section{Simulations}\label{simus}
In this section, we present a few simulations to understand more intuitively the structure of the dual chaos as introduced in section \ref{sec:KPZ}, i.e. $\gamma^2<2d$, $\gamma \overline{\gamma} =4$ and $\alpha=\frac{ \gamma^2}{4}$.

\begin{figure}[h]
\begin{center}
\includegraphics[width=15cm,height=6cm]{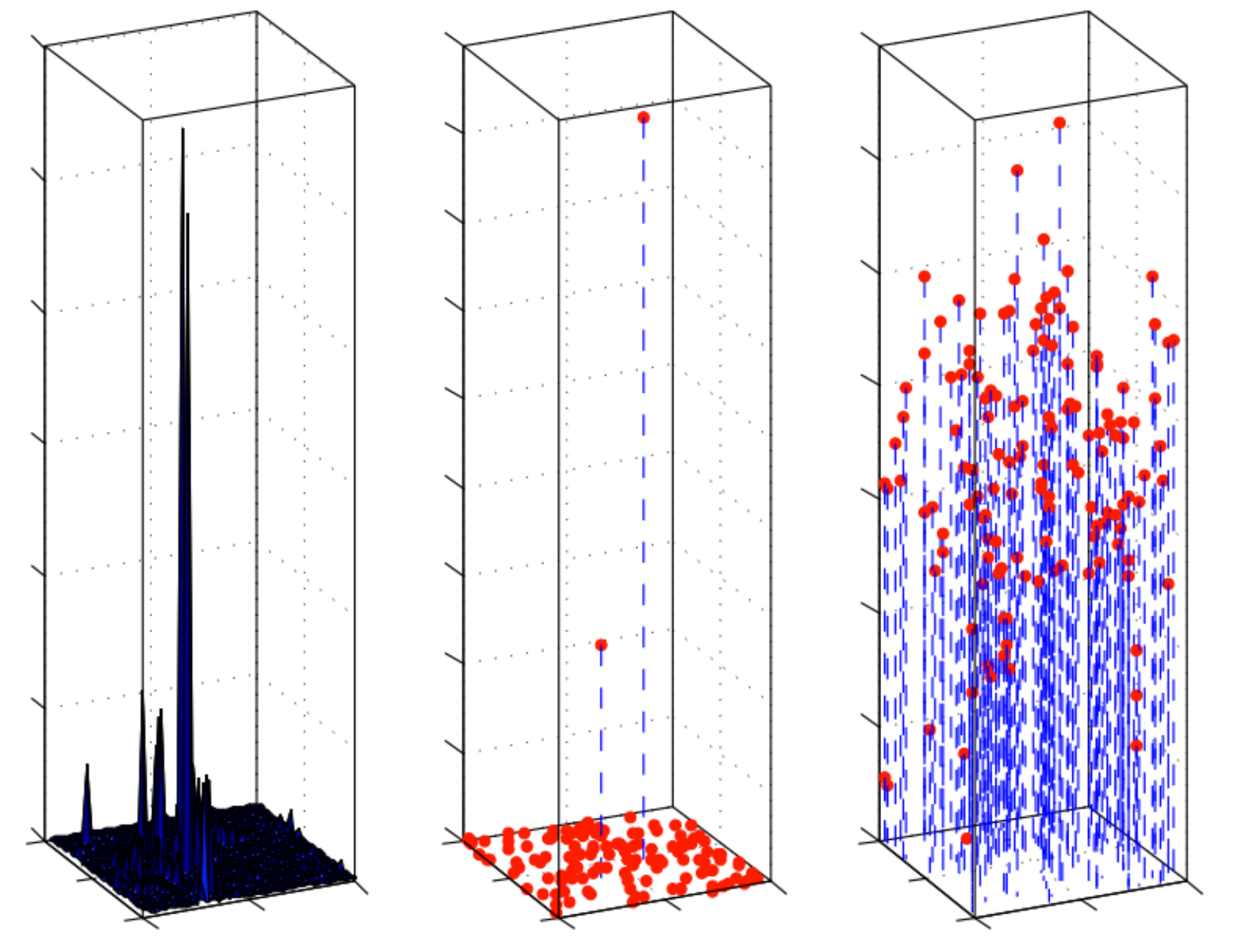}
\caption{Chaos and dual chaos for the value $\gamma^2=1$ (and then $\alpha=0.25$)}
\label{fig1}
\end{center}
\end{figure}

In Figure \ref{fig1}, we plot on the left hand side the ``density" of the usual chaos. The two other figures (middle and right) are concerned with the corresponding dual measures. In the middle, we plot the position and weights of the atoms of the dual measure. Notice that there are only a very small quantity of atoms with a very big weight. The other atoms have much smaller weights. To have a better picture of the values of these weights, we plot on the right-hand side the same picture with a logarithmic ordinate scale. 

\begin{figure}[h]
\centering
\subfloat[$\gamma^2=0.01$]{\includegraphics[width=0.47\linewidth]{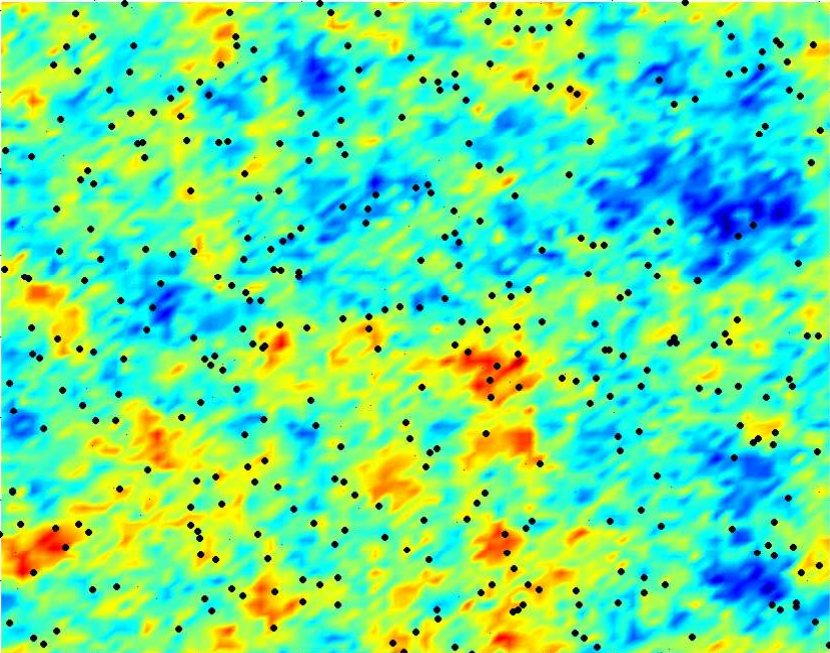}}
\,\,\subfloat[$\gamma^2=1$]{\includegraphics[width=0.47\linewidth]{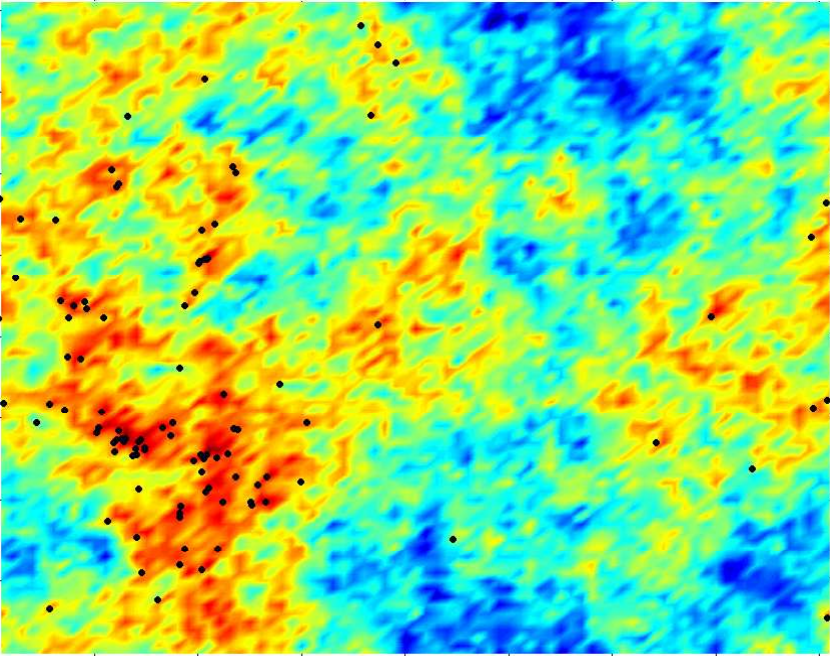}}\,\,
\subfloat[$\gamma^2=3.6$]{\includegraphics[width=0.47\linewidth]{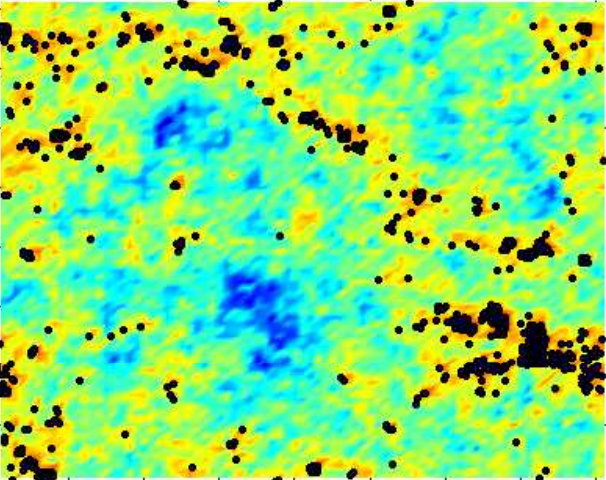}}
\caption{Spatial localization of a few big atoms  for different values of $\gamma$.}
\label{fig2}
\end{figure}

In Figure \ref{fig2}, we illustrate the influence of $\gamma$ on the spatial localization of the atoms. For different value of $\gamma$, we simulate a few atoms of the measure (the biggest). The colored background stands for the height profile of the measure $M$ plotted with a logarithmic intensity scale: red for areas with large mass and blue for areas with small mass. Localization of atoms is plotted in black. The larger $\gamma$ is, the more localized on  areas with large potential the atoms are.

\section{Perspectives}\label{pers}
 Here we develop a few comments and open problems related to this work.
 
\subsection{Dual chaos and possible renormalizations of degenerate Gaussian multiplicative chaos}
We continue to assume to be under exact scale invariance. For $\theta\ge 0$ consider  the associated sequence of measures 
$$
M_{\theta,n}(A)=\int_A e^{\theta X^n_x-\frac{\theta^2}{2}\mathbb{E}((X^n_x)^2)}\, dx,\quad n\ge 1,\ A\in B(\R^d). 
$$
Also define 
$$
\xi_\theta(q)=(d+\frac{\theta^2}{2})q-\frac{\theta^2}{2}q^2, \quad q\ge 0.  
$$
Recall that $(M_{\theta,n})_{n\ge 1}$ converges almost surely in the weak-star topology to a Radon measure $M_\theta$, which is almost surely positive or null according to whether $\theta^2<2d$ or $\theta^2\ge 2d$. 

By analogy with the study of Mandelbrot cascades and the fixed points of the associated smoothing transformation \cite{durrett,liu2000}, we may conjecture that when $\theta^2=2d$, the signed measures $-\frac{d}{d\theta} M_{\theta,n|\theta =\sqrt{2d}}$ weakly converge to a non-degenerate positive measure  $\widetilde M_{\sqrt{2d}}$. 

If $\theta^2>2d$, we have $\xi_\theta' (1)<0$ so that there exists a unique $\alpha\in (0,1)$ and a unique $\widetilde \alpha\in (\alpha,1)$ such that 
$$
\xi_\theta (\alpha)=d\quad\text{and}\quad \xi_\theta^*(\xi_\theta'(\widetilde \alpha))=-d,
$$ 
where $\xi^*_\theta(s)=\inf_{q\ge 0}s q-\xi_\theta(q)$. Indeed, the concavity of $\xi_\theta$ and the fact that $\xi_\theta(0)=0<\xi_\theta (1)=d$ and $\xi_\theta'(1)<0$ yields the existence and uniqueness of $\alpha$, at which we necessarily have $\xi_\theta'(\alpha)>\xi_\theta(\alpha)/\alpha=d/\alpha$. Then, we have $\xi_\theta^*(\xi_\theta'(\alpha))=\xi'( \alpha)\alpha-d>0$. Since $\xi_\theta^*(\xi_\theta'(1))=  \xi_\theta'(1)-d<-d$ and $\xi_\theta^*$ is concave, we get the existence and uniqueness of $\widetilde \alpha$. 

\medskip

Calculations  show that $\alpha=\frac{2d}{\theta^2}$ and  $\widetilde \alpha=\sqrt{\alpha}$.  Consequently, if we set $\overline\gamma=\theta$ and $\gamma= \alpha\overline\gamma$, we see that $\gamma\overline\gamma=2d$ and $\alpha$ is exactly the exponent used in the previous sections to establish the duality formula starting from the measure $M_\gamma$. Moreover, continuing the analogy with Mandelbrot cascades, the dual chaos $\overline M_\gamma$  is the expected non trivial solution, in ``replacement'' of $M_{\overline\gamma}$ which vanishes, of the equation 
\begin{equation}
\forall 0<\lambda<1,\quad (\overline{M}_\gamma(\lambda A))_{A\subset B(0,R)}\stackrel{law}{=}e^{\Omega'_\lambda}(\overline{M}(  A))_{A\subset B(0,R)}
\end{equation}
where $\Omega'_\lambda$ is a Gaussian random variable independent of the measure $(\overline{M}_\gamma(  A))_{A\subset B(0,R)}$ the law of which is characterized by:
$$\E[e^{q \Omega_\lambda}]=\lambda^{\xi_{\overline\gamma}(q)}$$
(for Mandelbrot cascades, in dimension 1, such measures have been identified as stable L\'evy subordinators in Mandelbrot time in \cite{BS}). 
\medskip

Thus, from the exact scale invariance point of view, the dual chaos of $M_\gamma$ provides a first way to renormalize  $M_{\overline\gamma}$, by giving a non trivial solution to the functional equation that would satisfy $M_{\overline\gamma}$ if it was not degenerate. 

Another way to build a non degenerate object from $M_{\theta,n}$ when $\theta^2> 2d$ is to consider the sequence of normalized measures $M_{\theta,n}/\|M_{\theta,n}\|$ (the equilibrium Gibbs measures considered for instance in \cite{CarDou}). Examining the behavior of $\|M_{\theta,n}\|$ shows that it approximately goes to 0 like $\exp (-n\xi_\theta'(\widetilde \alpha))$ when $\E[(X^n_0)^2]= n$. Then,  inspired by the recent progress made in this direction in the context of Branching random walks \cite{AidShi,BRV,madaule,Webb}, a tempting conjecture is that $M_{\theta,n}/\|M_{\theta,n}\|$ converges weakly, in law, to $\widetilde M^{(\widetilde \alpha)}_{\sqrt{2d}}/\|\widetilde M^{(\widetilde \alpha)}_{\sqrt{2d}}\|$, where $\widetilde M^{(\widetilde \alpha)}_{\sqrt{2d}}$ is defined  in the same way  as $\overline M$ was in~\eqref{defMbar}: fix a random measure $N_{\widetilde M_{\sqrt{2d}}}$  distributed on $\R^d\times\R^*_+$, and whose law conditionnally on $\widetilde M$ is that of a Poisson random measure with intensity $\frac{\widetilde M(dx)\,dz}{z^{1+\widetilde \alpha}}. 
$
Then,
\begin{align*} 
\forall A\in \mathcal{B}(\R^d),\quad \widetilde M^{(\widetilde \alpha)}_{\sqrt{2d}}(A)=&\int_A\int_{\R_+}z\,N_ {\widetilde M_{\sqrt{2d}}}(dx,dz). 
\end{align*}
Notice, however, that this second family of random measures cannot satisfy the duality relation since their power law spectrum, say $\xi_{\widetilde \alpha}$, satisfies $\xi_{\widetilde \alpha} (1)<d $. 
 
\begin{remark}
{\rm {\bf To continue the interpretation of duality as renormalization:} Consider an exact scale invariant log-infinitely divisible random measure (see \cite{bacry,rhovar} for precise definitions) for which we have enough exponential moments to discuss. We have a paramaterized family of multiplicative chaos $Q_{\theta,\epsilon}(r)=\exp (\theta\Lambda (C_r(t))/\E(\exp (\theta\Lambda (C_r(t)))$ and the associated measures $M_{\theta, \epsilon}$.  

Let $\psi(q)=\lim_{r\to 0} \log (\mathbb{E}( Q^q_{1,r}(x)))/\log (r)$. Now we have , 
$$
\xi_\theta(q)=(d+\psi(\theta))q-\psi(\theta q)
$$
and $\xi_\theta'(1)=d+\psi(\theta)-\theta\psi'(\theta)$. There is at most one positive (and at most one negative) solution $\theta_0$ to the equation $\xi_\theta'(1)=0$, and if $\theta_0>0$ exists, then for $\theta>0$ we have $\xi'_\theta(1)>0$ iff $\theta<\theta_0$. 

If $\overline \gamma=\theta>\theta_0$, we can have exactly the same discussion as in the Gaussian case, by considering the unique root $\alpha$ of $\xi_\theta(\alpha)=0$ and defining the dual chaos associated with  $\gamma=\alpha\theta$.  The dual KPZ relation is then naturally expressed via $\dim_{\overline M_\gamma} (K)= \alpha \dim_{M_\gamma} (K)$, of which the Gaussian case is a special case.

}
\end{remark}

\subsection{Singularity spectrum}
It would be interesting to compute the free energy of the measure $\overline{M}$, namely proving that the following limit is not trivial:
$$\lim_{n\to \infty}2^{n(\overline{\xi}(q)-d)}\sum_{I\in C_n}\overline{M}(I)^q\quad (q\in \R),$$
where $C_n$ stands for the set of all dyadic cubes (included in the unit cube) with side length $2^{-n}$. This thermodynamic point of view is closely related to the calculation of the $L^q$-spectrum of the measure $\overline M$, defined as  
$$
q\in\mathbb{R}\mapsto \tau_{\overline M}(q)=\liminf_{r\to 0^+}\frac{\log \sup\Big \{\sum_i \overline M(B(x_i,r))^q\Big \}}{\log(r)},
$$
where the supremum is taken over all the centered packing of $[0,1]^d$ by closed balls of radius  $r$. By analogy with the study achieved in \cite{BS}, we conjecture that on the one hand,
$$
\tau_{\overline M}(q)= 
\begin{cases}
\overline \xi'(q_-) q&\text{if } q\le q_-,\\
\overline \xi(q)-d&\text{if } q_-\le q\le \alpha,\\
0&\text{if }q\ge \alpha,
\end{cases}
$$
where $q^-$ is the unique negative solution of $\overline \xi^*(\overline\xi'(q))=-d$, and on the other hand that the multifractal formalism holds for $\overline M$: defining $$E_\delta=\big\{x\in [0,1]^d;\liminf_{r\to 0^+}\frac{\ln \overline{M}(B(x,r))}{\ln (r)}=\delta\big\}
\quad (\delta\ge 0),$$
with probability 1, the singularity spectrum of $\overline M$, i.e. the mapping $\delta\ge 0\mapsto \dim E_\delta$, is given by $\delta\ge 0\mapsto \tau_{\overline M}^*(\delta)=\inf\{\delta q-\tau_{\overline M}(q):q\in \R\}$, a negative dimension meaning that $E_\delta=\emptyset$. 

\appendix


\section{Examples of sigma-positive kernels}\label{rem:pas}
In this section, we detail a few examples of sigma-positive kernels, apart from  the Green function already explained in subsection \ref{sub:GFF}. More precisely, we give two different classes of sigma-positive kernels, which yield two different notions of stochastic scale invariance for the associated Gaussian multiplicative chaos.

\subsection{Exact stochastic scale invariance}

In this section, we describe how to construct kernels yielding the exact scale invariance relations of Theorem \ref{perfect}. This is useful in computations and it is possible to deduce all the other situations from this one.

We define on $\R_+$ the measure $\nu_T(dt)=\ind_{[0,T]}(t)\frac{dt}{t^2}+\frac{1}{T}\delta_T(dt)$ where $\delta_x$ denotes the Dirac mass at $x$. For $\mu>0$, it is straightforward to check that
\begin{equation}\label{formln}
\forall x\in\R^d,\quad \ln_+\frac{T}{|x|} =\frac{1}{\mu}\int_0^{+\infty}(t-|x|^\mu)_+\nu_{T^\mu}(dt).
\end{equation}

-In dimension $d=1$, it is straightforward to check that the function $x\mapsto (t-|x|)_+$ is of positive type. So, the kernel $K(x)=\gamma^2\ln_+\frac{T}{|x|}$ is of sigma positive type. The kernels $k_n$ can be easily computed:
$$k_n(x)=\left\{\begin{array}{ll}
0 &\text{if }|x|>T,\\
\gamma^2\ln_+\frac{T}{|x|} & \text{if }\frac{T}{n}\leq |x|\leq T,\\
 \gamma^2\ln n+ \big(1- \frac{n|x|}{T}\big)                &\text{if } 0\leq |x|\leq \frac{T}{n}.
 \end{array}\right.$$

-In dimension $d=2$, Pasenchenko \cite{cf:pas} proved that the function $(1-|x|^{1/2})_+ $ is positive definite in dimension $2$. Choosing $\mu=2$ in \eqref{formln}, we can thus write 
$$\forall x\in\R^2,\quad \gamma^2\ln_+\frac{T}{|x|} =\sum_{n\geq 1} q_n(x),$$ where $q_n$ is the continuous positive and positive definite kernel $$\forall x\in\R^2,\quad q_n(x)=2\gamma^2\int_{\frac{T^{1/2}}{n^{1/2}}}^{\frac{T^{1/2}}{(n-1)^{1/2}}}(t-|x|^\mu)_+\nu_{T^{1/2}}(dt).$$
A simple computation shows that
$$k_n(x)=\left\{\begin{array}{ll}
0 &\text{if }|x|>T,\\
\gamma^2\ln_+\frac{T}{|x|} & \text{if }\frac{T}{n}\leq |x|\leq T,\\
 \gamma^2\ln n+2\big(1-\sqrt{\frac{n|x|}{T}}\big)                &\text{if } 0\leq |x|\leq \frac{T}{n}.
 \end{array}\right.$$
 
 -In dimension $d\geq 3$, it is proved in \cite{rhovar} that there exists a continuous bounded function $g:\R^d\to \R$, constant in a neighborhood of $0$ such that
\begin{equation}\label{perfkernel}
K(x)=\gamma^2\ln_+\frac{T}{|x|} +g(x)
\end{equation}
 is of sigma positive type.

\subsection{$\star$-scale invariance}
A simple way of constructing sigma positive kernels is given by
\begin{equation}\label{star}
\forall x\in\R^d,\quad K(x)=\int_1^{\infty}\frac{k(xu)}{u}\,du,
\end{equation}
 where $k$ is a continuous positive kernel of positive type. Such kernel  is of sigma positive type since the decomposition can be realized by
  $$q_n(x)=\int_{2^n}^{2^{n+1}}\frac{k(xu)}{u}\,du.$$
  Furthermore, $K$ takes on the form \eqref{formgam} with $\gamma^2=k(0)$. Such kernels are related to the notion of $\star$-scale invariance (see \cite{allez,sohier}).
\section{Proofs of Section \ref{sec:chaos}}
\subsection*{Preliminary computations}
In this section, we prove a preliminary result, which is about the convergence in law of the sequence $(\overline{M}_n)_n$. We will use the following relation valid for any $0<\beta<1$ and $x\geq 0$:
\begin{equation}\label{fourier}
x^\beta= \frac{\beta}{\Gamma(1- \beta)}\int_0^{\infty}(1-e^{-xz})\frac{dz}{z^{1+\beta}}.
\end{equation}
Therefore, we have for all $u\geq 0$:
\begin{equation}\label{expn}
 \E[e^{-u \overline{M}_n(A)}]=\E\Big[e^{\int_A\int_{\R_+}(e^{-zu}-1)\frac{1}{z^{1+\alpha}}\,dz M_n(dx)}\Big]=\E[e^{-\frac{\Gamma(1-\alpha)}{\alpha}u^\alpha M_n (A)}].
\end{equation}
Similarly, we have
\begin{equation}\label{expmultin}
\E[e^{-u_1 \overline{M}_n(  A_1)+\dots-u_p\overline{M}_n(  A_p)}]=\E\big[e^{-\frac{\Gamma(1-\alpha)}{\alpha}\big(u_1^\alpha M_n(A_1)+\dots+u_p^\alpha M_n(A_p)\big)}\big]
\end{equation} valid for all $u_1,\dots,u_p\in\R_+$ and all disjoint Borelian subsets $A_1,\dots,A_p\subset \R^d$. This is just a consequence of \eqref{expn} and the fact that the random measure $n_\alpha$ is independently scattered (see the theory of independently scattered random measures in  \cite{Ros}). 

Since $(M_n)_n$ almost surely weakly converges towards $M$, we get as a direct consequence  the convergence in law of the finite marginals of  $((\overline{M}_n(A))_{A\in \mathcal B(\R^d)})_n$. Since any compactly supported continuous function can be uniformly approximated from above and below by linear combinations of indicator functions of disjoint semi-open dyadic cubes, we get the convergence in law of the finite marginals of $(\overline{M}_n(f))_{f\in C_c(\R^d)}$. Moreover, for each compact cube $A_N=[-N,N]^d$, $N\ge 1$,  the sequence $(\overline{M}_n(A_N)_n$ is tight, so the laws of the  Borel measures $\overline{M}_n$ restricted to $A_N$ form a tight sequence, which  implies the tightness of the laws of the Radon measures $\overline{M}_n$. Thus  we obtain that the sequence $(\overline{M}_n)_n$ converges in law towards a random Radon measure, call it $\overline{M}$, as well as a characterization of the law of the limiting measure:
\begin{equation}\label{expmulti}
\E[e^{-u_1 \overline{M}(  A_1)-\dots-u_n\overline{M}(  A_p)}]=\E\big[e^{-\frac{\Gamma(1-\alpha)}{\alpha}\big(u_1^\alpha M(A_1)+\dots+u_n^\alpha M(A_n)\big)}\big]
\end{equation} valid for all $u_1,\dots,u_p\in\R_+$ and all disjoint Borelian subsets $A_1,\dots,A_p\subset \R^d$. In particular, we have for all $u\geq 0$:
\begin{equation}\label{exp}
 \E[e^{-u \overline{M}(A)}]=\E[e^{-\frac{\Gamma(1-\alpha)}{\alpha}u^\alpha M (A)}].
\end{equation}
This proves the relations \eqref{relfundn}. Let us prove  \eqref{relmom}. For $0<\beta<\alpha$, we have
\begin{align*}
\E[(\overline{M}(A))^\beta]=&\frac{\beta}{\Gamma(1- \beta)}\int_0^{\infty}\big(1-\E[e^{-w\overline{M}(A)}]\big)\frac{dw}{w^{1+\beta}}\\
=& \frac{\beta}{\Gamma(1- \beta)} \int_0^{\infty}\big(1-\E[e^{-\frac{\Gamma(1-\alpha)}{\alpha} w^\alpha M(A)}]\big)\frac{dw}{w^{1+\beta}}.
 \end{align*}
We make the change of variables $y=w^\alpha$ to get:
\begin{align}\label{momn}
\E[(\overline{M}(A))^\beta]=& \frac{\beta}{\alpha\Gamma(1- \beta)} \int_0^{\infty}\big(1-\E[e^{-\frac{\Gamma(1-\alpha)}{\alpha}y M(A)}]\big)\frac{dy}{y^{1+\frac{\beta}{\alpha}}}\\
=&  \frac{  \Gamma(1-\beta/\alpha)\Gamma(1-\alpha)^{\beta/\alpha}}{ \Gamma(1- \beta)\alpha^{\beta/\alpha}}  \E[(M(A))^\frac{\beta}{\alpha}].\nonumber
 \end{align}

Notice that the measure defined by \eqref{defMbar} obviously satisfies \eqref{relfundn}. Indeed, let us consider the measure $\overline{M}$ defined by \eqref{defMbar}. Let us also consider $u_1,\dots,u_p\in\R_+$ and disjoint Borelian subsets $A_1,\dots,A_p\subset \R^d$.  Observe that, conditionally on $M$, the random variables $\overline{M} (  A_1),\dots,\overline{M} (  A_p)$ are independent. Thus, by using \eqref{fourier},  we get
\begin{align*}
 \E[e^{-u_1 \overline{M} (  A_1)-\cdots-u_p\overline{M} (  A_p)}]&=\E\big[\E[e^{-u_1 \overline{M} (  A_1)-\cdots-u_p\overline{M} (  A_p)}|M]\big]\\
 &=\E\big[e^{-\frac{\Gamma(1-\alpha)}{\alpha}\big(u_1^\alpha M(A_1)+\dots+u_p^\alpha M(A_p)\big)}\big].
\end{align*}
Therefore the construction of Theorem \ref{as}  yields the same measure in law as the measure defined by \eqref{defMbar}. 

\vspace{3mm}
\noindent {\it Proof of Proposition \ref{def}.} From \eqref{relfundn}, we deduce that the law of $\overline{M}$ is characterized by that of $M$, which does not depend on the chosen decomposition (see \cite{Kah}). Furthermore, since $N_M$ (define before \eqref{defMbar}) is a Poisson random measure conditionally to $M$, it is clear that $\overline{M}$ (see \eqref{defMbar}) is almost surely purely atomic. \qed

\subsection*{Proofs of Theorem \ref{as}}
 
 Now we tackle the convergence in probability. We will prove that we can extract from each subsequence a subsequence converging in probability. 
 
 \medskip
 
\noindent (1) Consider a subsequence $(\phi(n))_n$. Since ${\rm Var}(X^{\phi(n)}_0)=\sum_{k=1}^{\phi(n)}{\rm Var}(Y^k_0)\to \infty$ as $n\to \infty$, we can find a subsequence $(\psi(n))_n$ such that
\begin{equation}\label{var}
\sum_{k=\psi(n-1)+1}^{\psi(n)}{\rm Var}(Y^k_0)\geq \rho\ln n
 \end{equation}
  where $\rho>\frac{8}{\gamma^2}$. Let us  now prove that the sequence $(\overline{M}_{\psi(n)}(A))_n$ converges in probability for every bounded Borelian subset $A$ of $\R^d$. 
 
We can rearrange the $(Y^n)_n$ and set 
\begin{equation}\label{prime}
Y^{'n}=\sum_{k=\psi(n-1)+1}^{\psi(n)} Y^k,\quad X^{'n}=\sum_{k=1}^{n} Y^{'k}.
\end{equation}
For the sake of clarity, we will omit the superscript $'$ from the notations. So we will just assume below that the processes $(X^n)_n,(Y^n)_n$ satisfy the usual properties of Section \ref{back} together with  the constraint ${\rm Var}(Y^n_0)\geq \rho\ln n$ .

We fix $c>0$. We denote by $\mathcal{F}_n$ ($n\geq 0$) the sigma algebra generated by the random measure $n_\alpha$ and the random processes $(X^{p})_{p\leq  n}$.   We set
\begin{align}
\overline{M}_n^c(A)=&\int_A\int_0^{c e^{\frac{\gamma}{\alpha} X^{n }_x-\frac{\gamma^2}{2\alpha}\E[(X^{n }_x)^2]}}e^{\frac{\gamma}{\alpha} X^{n }_x-\frac{\gamma^2}{2\alpha}\E[(X^{n }_x)^2]}z\,n_\alpha(dx,dz)\label{trunc}\\
=&\int_A\int_0^{+\infty}\ind_{\{0\leq z \leq c e^{\frac{\gamma}{\alpha} X^{n }_x-\frac{\gamma^2}{2\alpha}\E[(X^{n }_x)^2]}\}}e^{\frac{\gamma}{\alpha} X^{n }_x-\frac{\gamma^2}{2\alpha}\E[(X^{n }_x)^2]}z\,n_\alpha(dx,dz)\nonumber
\end{align}
We have 
\begin{align*}
\E\Big[&(\overline{M}_{n+1}^c(A))^\alpha  |\mathcal{F}_n\Big]\\ 
= &\E\Big[\Big(\int_A\int_0^{c e^{\frac{\gamma}{\alpha} X^{n+1 }_x-\frac{\gamma^2}{2\alpha}\E[(X^{n +1}_x)^2]}}e^{\frac{\gamma}{\alpha} X^{n+1 }_x-\frac{\gamma^2}{2\alpha}\E[(X^{n+1 }_x)^2]}z\,n_\alpha(dx,dz)\Big)^\alpha|\mathcal{F}_n\Big]\\
\geq &\E\Big[\Big(\int_A\int_0^{c e^{\frac{\gamma}{\alpha} X^{n+1 }_x-\frac{\gamma^2}{2\alpha}\E[(X^{n +1}_x)^2]}}\ind_{\big\{\frac{\gamma}{\alpha} Y^{n+1 }_x-\frac{\gamma^2}{2\alpha}\E[(Y^{n+1 }_x)^2]\geq 0 \big\}}e^{\frac{\gamma}{\alpha} X^{n+1 }_x-\frac{\gamma^2}{2\alpha}\E[(X^{n+1 }_x)^2]}z\,n_\alpha(dx,dz)\Big)^\alpha|\mathcal{F}_n\Big]\\
\geq &\E\Big[\Big(\int_A\int_0^{c e^{\frac{\gamma}{\alpha} X^{n }_x-\frac{\gamma^2}{2\alpha}\E[(X^{n  }_x)^2]}}\ind_{\big\{\frac{\gamma}{\alpha} Y^{n+1 }_x-\frac{\gamma^2}{2\alpha}\E[(Y^{n+1 }_x)^2]\geq 0 \big\}}e^{\frac{\gamma}{\alpha} X^{n+1 }_x-\frac{\gamma^2}{2\alpha}\E[(X^{n+1 }_x)^2]}z\,n_\alpha(dx,dz)\Big)^\alpha|\mathcal{F}_n\Big],
\end{align*}
using the fact that $\frac{\gamma}{\alpha} X^{n+1 }_x-\frac{\gamma^2}{2\alpha}\E[(X^{n +1}_x)^2]=\frac{\gamma}{\alpha} Y^{n+1 }_x-\frac{\gamma^2}{2\alpha}\E[(Y^{n+1 }_x)^2]+ \frac{\gamma}{\alpha} X^{n}_x-\frac{\gamma^2}{2\alpha}\E[(X^{n}_x)^2]$. Since the mapping $x\in\R_+ \mapsto x^\alpha$ is concave, we apply Jensen's inequality with respect to the restriction to $A$ of the  measure $\overline{M}_{n}^c/\overline{M}_{n}^c(A)$, conditionally on  $\overline{M}_{n}^c(A)\neq 0$, and  get:
\begin{align*}
\E\Big[&(\overline{M}_{n+1}^c(A))^\alpha  |\mathcal{F}_n\Big]\\ 
\geq &\int_A\int_0^{c e^{\frac{\gamma}{\alpha} X^{n }_x-\frac{\gamma^2}{2\alpha}\E[(X^{n  }_x)^2]}}\E\Big[e^{\gamma Y^{n+1 }_x-\frac{\gamma^2}{2 }\E[(Y^{n+1 }_x)^2]}\ind_{\big\{\frac{\gamma}{\alpha} Y^{n+1 }_x-\frac{\gamma^2}{2\alpha}\E[(Y^{n+1 }_x)^2]\geq 0 \big\}}\Big]\times\dots\\
&\dots\times e^{\frac{\gamma}{\alpha} X^{n }_x-\frac{\gamma^2}{2\alpha}\E[(X^{n }_x)^2]}z\,n_\alpha(dx,dz) \,\,\ind_{\{\overline{M}_{n}^c(A)\neq 0\}}(\overline{M}_{n }^c(A))^{\alpha-1}\\
=&\ind_{\{\overline{M}_{n}^c(A)\neq 0\}}(\overline{M}_n^c(A))^\alpha \E\Big[e^{\gamma Y^{n+1 }_0-\frac{\gamma^2}{2 }\E[(Y^{n+1 }_0)^2]}\ind_{\big \{\frac{\gamma}{\alpha} Y^{n+1 }_0-\frac{\gamma^2}{2\alpha}\E[(Y^{n+1 }_0)^2]\geq 0 \big\}}\Big]\\
=& (\overline{M}_n^c(A))^\alpha \E\Big[e^{\gamma Y^{n+1 }_0-\frac{\gamma^2}{2 }\E[(Y^{n+1 }_0)^2]}\ind_{\big\{\frac{\gamma}{\alpha} Y^{n+1 }_0-\frac{\gamma^2}{2\alpha}\E[(Y^{n+1 }_0)^2]\geq 0 \big\}}\Big].
\end{align*}
By using a Girsanov transform we deduce
\begin{align*}
\E\Big[(\overline{M}_{n+1}^c(A))^\alpha  |\mathcal{F}_n\Big]\geq & (\overline{M}_n^c(A))^\alpha\,\,\P\Big(  \gamma  Y^{n+1 }_0+\frac{\gamma^2}{2 }\E[(Y^{n+1 }_0)^2]\geq 0\Big)\\
=&(\overline{M}_n^c(A))^\alpha\,\,\Big(1-\P\big(  \mathcal{N}(0,1)\geq\frac{\gamma}{2 }{\rm Var}(Y^{n+1 }_0)^{\frac{1}{2}}\big)\Big).
\end{align*}
Let us set $b_n=\P\big(  \mathcal{N}(0,1)\geq\frac{\gamma}{2 }{\rm Var}(Y^{n+1 }_0)^{1/2}\big)$ in such a way that 
\begin{align*}
\E\Big[(\overline{M}_{n+1}^c(A))^\alpha  |\mathcal{F}_n\Big]\geq  (\overline{M}_n^c(A))^\alpha\,\,\Big(1-b_n\Big).
\end{align*}
We further stress that \eqref{var} implies that the series $\sum_nb_n$ is absolutely convergent.
Let us define 
$$c_n=\prod_{k=1}^n\frac{1}{1-b_k}$$
and $$A_n=c_n (\overline{M}_n^c(A))^\alpha.$$
The sequence $(A_n)_n$ is a positive sub-martingale. Let us prove that it is bounded in $L^1$. We have
\begin{align*}
\E[A_n]=&c_n\E\Big[(\overline{M}_{n }^c(A))^\alpha \Big]\\=&\frac{\alpha c_n}{\Gamma(1-\alpha)}\int_0^{+\infty}\frac{1-\E[e^{-u \overline{M}_{n }^c(A)}]}{u^{1+\alpha}}\,du\\
=& \frac{\alpha c_n}{\Gamma(1-\alpha)}\int_0^{+\infty}\Big(1-\E\Big[e^{\int_A\int_0^{c e^{\frac{\gamma}{\alpha} X^{n }_x-\frac{\gamma^2}{2\alpha}\E[(X^{n }_x)^2]}}  \big(e^{-uze^{\frac{\gamma}{\alpha} X^{n }_x-\frac{\gamma^2}{2\alpha}\E[(X^{n }_x)^2]}}-1\big) \frac{1}{z^{1+\alpha}}\,dz   dx }\Big]\Big)\frac{1}{u^{1+\alpha}}\,du\\
=& \frac{\alpha c_n}{\Gamma(1-\alpha)}\int_0^{+\infty}\Big(1-\E\Big[e^{\int_A\int_0^{c }  \big(e^{-uy}-1\big) \frac{1}{y^{1+\alpha}}\,dy  \, e^{\gamma X^{n }_x-\frac{\gamma^2}{2 }\E[(X^{n }_x)^2]}dx }\Big]\Big)\frac{1}{u^{1+\alpha}}\,du\\
=& \frac{\alpha c_n}{\Gamma(1-\alpha)}\int_0^{+\infty}\Big(1-\E\Big[e^{-M_n(A) L_c(u)}\Big]\Big)\frac{1}{u^{1+\alpha}}\,du,
\end{align*}
 where $L_c(u)=\int_0^{c }  \big(1-e^{-uy}\big) \frac{1}{y^{1+\alpha}}\,dy$. Notice that $L_c(u)\geq 0$ for each $u\geq 0$. Since the mapping $x\mapsto -e^{-x L_c(u)}$ is concave, we can use Jensen's inequality to get
\begin{align*}
\E[A_n]\leq & \frac{\alpha c_n}{\Gamma(1-\alpha)}\int_0^{+\infty}\Big(1-\E\Big[e^{-|A|L_c(u)}\Big]\Big)\frac{1}{u^{1+\alpha}}\,du\\
=&c_n\E\Big[\Big(\int_A\int_0^cz\,n_\alpha(dx,dz)\Big)^\alpha\Big].
\end{align*}
The expectation in the above right-hand side is finite. Furthermore, the convergence of the series $\sum_nb_n$ implies the convergence of the sequence $(c_n)_n$ towards $\prod_{k=1}^{+\infty}\frac{1}{1-b_n}\in]0,+\infty[$. Therefore, the sub-martingale $(A_n)_n$ almost surely converges. So does $(\overline{M}_{n }^c(A))_n$. Let us denote by $\overline{M}^c(A)$ its limit. Obviously, the mapping $c\mapsto \overline{M}^c(A)$ is increasing and thus converges as $c$ goes to $\infty$. Let us denote by $\overline{M}^\infty(A)$ the limit.

Now we prove that the sequence $(\overline{M}_n(A))_n$  converges in probability towards $\overline{M}^\infty(A)$. We have for $\delta>0$:
\begin{align}
 \P(&|\overline{M}_n(A)-\overline{M}^\infty(A)|>\delta)\nonumber\\&\leq  \P(|\overline{M}_n(A)-\overline{M}^c_n(A)|>\frac{\delta}{3})+ \P(|\overline{M}^c_n(A)-\overline{M}^c(A)|>\frac{\delta}{3})+ \P(|\overline{M}^c(A)-\overline{M}^\infty(A)|>\frac{\delta}{3})\nonumber\\
 &\leq  \frac{3^\beta}{\delta^\beta}\E[|\overline{M}_n(A)-\overline{M}^c_n(A)|^\beta]+ \P(|\overline{M}^c_n(A)-\overline{M}^c(A)|>\frac{\delta}{3})+ \P(|\overline{M}^c(A)-\overline{M}^\infty(A)|>\frac{\delta}{3})\label{prob:de}
\end{align}
where $\beta<\alpha$. We evaluate the first quantity. We have:
\begin{align*}
\E[|\overline{M}_n(A)&-\overline{M}^c_n(A)|^\beta]\\
=&\frac{\beta}{\Gamma(1-\beta)}\int_0^{+\infty}\frac{1-\E[e^{-u (\overline{M}_n(A)-\overline{M}^c_n(A))}]}{u^{1+\beta}}\,du\\
=&\frac{\beta}{\Gamma(1-\beta)}\int_0^{+\infty}\Big(1-\E\Big[e^{\int_A\int^{+\infty}_{c e^{\frac{\gamma}{\alpha} X^{n }_x-\frac{\gamma^2}{2\alpha}\E[(X^{n }_x)^2]}}  \big(e^{-uze^{\frac{\gamma}{\alpha} X^{n }_x-\frac{\gamma^2}{2\alpha}\E[(X^{n }_x)^2]}}-1\big) \frac{1}{z^{1+\alpha}}\,dz   dx }\Big]\Big)\frac{1}{u^{1+\beta}}\,du\\
=& \frac{\beta}{\Gamma(1-\beta)}\int_0^{+\infty}\Big(1-\E\Big[e^{\int_A\int_c^{+\infty }  \big(e^{-uy}-1\big) \frac{1}{y^{1+\alpha}}\,dy\,   e^{\gamma X^{n }_x-\frac{\gamma^2}{2 }\E[(X^{n }_x)^2]}dx }\Big]\Big)\frac{1}{u^{1+\beta}}\,du\\
=&\frac{\beta}{\Gamma(1-\beta)}\int_0^{+\infty}\Big(1-\E\Big[e^{-M_n(A) U_c(u)}\Big]\Big)\frac{1}{u^{1+\beta}}\,du
\end{align*}
where $U_c(u)=\int_c^{+\infty }  \big(1-e^{-uy}\big) \frac{1}{y^{1+\alpha}}\,dy$. Notice that $U_c(u)\geq 0$ for each $u\geq 0$. From Jensen's inequality again, we have
\begin{align}
\E[|\overline{M}_n(A)-\overline{M}^c_n(A)|^\beta]\leq & \frac{\beta}{\Gamma(1-\beta)}\int_0^{+\infty}\Big(1-\E\Big[e^{-|A|U_c(u)}\Big]\Big)\frac{1}{u^{1+\beta}}\,du\nonumber\\
=& \E\Big[\Big(\int_A\int_c^{+\infty}z\,n_\alpha(dx,dz)\Big)^\beta\Big].\label{cvun}
\end{align}
This latter quantity converges to $0$ as $c$ goes to $\infty$ (uniformly with respect to $n$).

Now we come back to \eqref{prob:de} to complete the proof.
We can fix $c>0$ so as to make the first and third quantity as small as we please. Indeed, concerning the first quantity, it results from the bound just above \eqref{cvun} and concerning the third quantity, it results from the almost sure convergence of $\overline{M}^c(A)$ towards $\overline{M}^\infty(A)$. For such a $c$, we can find $N$ such that the second quantity is also as small as we please for $n\geq N$. This yields the convergence in probability of $(M_{\psi(n)}(A))_n$. 

Now suppose that $(\overline M_{\psi'(n)}(A))_n$ is another subsequence of $(\overline M_n(A))_n$ converging in probability. Build a new subsequence as follows: first pick a term in the sequence $(\psi(n))_n$, then pick a term in  $(\psi'(n))_n$ with large index, then pick a  large term in  $(\psi(n))_n$, and so on with gaps large enough so that the resulting sequence $(\overline M_{\xi(n)})_n$ satisfies condition \eqref{var}. It follows that $(\overline M_{\xi(n)})_n$ is extracted both from $(\overline M_{\psi(n)}(A))_n$ and $(\overline M_{\psi'(n)}(A))_n$, and it converges in probability. The limits in probability along  $(\overline M_{\psi(n)}(A))_n$ and $(\overline M_{\psi'(n)}(A))_n$ must thus be the same.

Observe that we have just proved the convergence in probability of the sequence $(\overline{M}_n(A))_n$ for each Borelian bounded set $A$. If $f$ is a compactly supported continuous function, it is straightforward to do the same for $(\overline M_n(f))_n=(\int f(x)\overline{M}_n(dx))_n$ (the computations and estimates are similar). For each $N\in \N$, let $(f_{N,j})_{j\in\N}$ be a dense sequence in the space of  continuous functions compactly  supported in $[N,N]^d$. Since all the sequences $(\overline M_n(f_{N,j}))_n$ converge in probability, using the Cantor diagonal extraction principle we can find a subsequence $(\overline M_{\psi(n)})_n$ such that, with probability 1, for all $N,j\in \N$, $(\overline M_{\psi(n)}(f_{N,j})_n$ converges to, say, $\overline M(f_{N,j})$. It is then standard that $\overline M$ can be extended to a random Radon measure on $\R^d$ because of its positivity. Still with the same arguments as above, we also have that $\big(\overline{M}_n(f_{N,j}))_{(N,j) \in\mathbb{N}^2}\big)_n$ converges in probability to $\big(\overline{M}(f_{N,j}))_{(N,j) \in\mathbb{N}^2}\big)_n$, where  the space  $\R^{\mathbb{N}^2}$ is endowed with the product topology. Now let us prove that $(\overline M_n)_n$ weakly converges in probability to $\overline M$. For this we endow the space of Radon measures on $\R^d$ with the standard distance $d(\mu,\mu')=\sum_{N\ge 0}2^{-N-1}\sum_{j\ge 0}2^{-j-1} \min(1, |\mu (f_{N,j})-\mu'(f_{N,j})|)$. Fix $\epsilon>0$ and $N_0\in\N_+$ such that $2^{-N_0}\le \epsilon/4$. Now fix $\eta>0$ and then $n_0\in\N$ such that for all $n\ge n_0$, for all $0\le N\le N_0$ and $0\le j\le N_0$ we have $\P (|\overline M_n(f_{N,j})-\overline M(f_{N,j})|\ge \epsilon/2)\le  (N_0+1)^{-2} \eta $. It is direct to see that $\{d(\overline M_n,\overline M)\ge \epsilon\}\subset \bigcup_{0\le N\le N_0,\, 0\le j\le N_0} \{ |\overline M_n(f_{N,j})-\overline M(f_{N,j})|\ge \epsilon/2\}$, whose probability is smaller than $\eta$. We thus showed the convergence in probability of $(\overline M_n)_n$, hence statement~(2).  

\smallskip

\smallskip
\noindent Statement (3) and (4) follow from the convergence in law of $(\overline M_n)_n$ to a limit characterized by  the relation  \eqref{expmulti}, the limiting law being non trivial if and only if $M$ is non trivial (the necessary and sufficient condition $\gamma^2<2d$ for non triviality of $M$ is proved in \cite{Kah}).
 \qed

\begin{remark}
The reader may find the above proof more tricky than expected. Actually, the truncation suggested in \eqref{trunc} is not the more natural way that we may think of to tackle the problem. The first idea that we may come up with is rather to define
\begin{align*}
\overline{M}_n^c(A)=&\int_A\int_0^{c }e^{\frac{\gamma}{\alpha} X^{n }_x-\frac{\gamma^2}{2\alpha}\E[(X^{n }_x)^2]}z\,n_\alpha(dx,dz) .
\end{align*}
It is straightforward to check that $(\overline{M}_n^c(A)^\alpha)_n$ is a submartingale. But it is not bounded in $L^1$. Indeed, if it was, its limit would be $\overline{M}(A)^\alpha$ regardless of the value of $c$ because the chaos ``kills" the big jumps, i.e. $\int_A\int_c^{\infty }e^{\frac{\gamma}{\alpha} X^{n }_x-\frac{\gamma^2}{2\alpha}\E[(X^{n }_x)^2]}z\,n_\alpha(dx,dz)\to 0$ as $n\to \infty$ almost surely. Thus, $\overline{M}(A)$ would admit a moment of order $\alpha$, which is impossible (see below).
\end{remark}

\subsection*{Proofs of Proposition \ref{prop:mom}}
For $\beta<\alpha$, we can use relation \eqref{momn} to show  the existence of the moments and the dual relation \eqref{relmom}. If $\overline{M}$ possesses a moment of order $\alpha$ then the left-hand side of equation \eqref{relmom} must converge as $\beta\to \alpha$. But it is equal to the right-hand side, which diverges because of the term $ \Gamma(1-\beta/\alpha)$ and the fact that the measure $M$ possesses a non trivial moment of order $1$.

\subsection*{Proof of Theorem \ref{perfect} }

 First we stress that it has already been proved that the chaos measure  $M$, associated to the given kernel $K$,  satisfies the scale invariance relation (see \cite{rhovar}) for some $R>0$:
\begin{equation*} 
\forall 0<\lambda<1,\quad (M(\lambda A))_{A\subset B(0,R)}\stackrel{law}{=}\lambda^{d }e^{ \Omega_\lambda }(M(  A))_{A\subset B(0,R)}
\end{equation*}
where $\Omega_\lambda$ is a Gaussian random variable independent of the measure $(\overline{M}(  A))_{A\subset B(0,R)}$ the law of which is characterized by:
$$\E[e^{q \Omega_\lambda}]=\lambda^{\frac{\gamma^2}{2}q-\frac{\gamma^2}{2}q^2}.$$ The results then easily follows from the relation
$$\E[e^{-u_1 \overline{M}(  A_1)+\dots-u_n\overline{M}(  A_n)}]=\E[e^{-u_1^\alpha M(A_1)+\dots-u_n^\alpha M(A_n)}]$$ valid for all $u_1,\dots,u_n\in\R$ and all disjoint Borelian subsets $A_1,\dots,A_n\subset \R^d$. 
\qed

\subsection*{Proof of Corollary \ref{notperf}}
 Let us write the kernel $K$ as
$$K(x)=K_p(x)+h(x)$$ where $K_p$ is the ``perfect kernel" given by \eqref{perfkernel} and $g$ is some continuous bounded function over $\R^d$. Even if it means adding to $K$ a constant, we may assume that $h(0)=0$ and, without loss of generality, we assume $R=1$.
For $t>0$, we define 
$$G_t=\sup_{|x|\leq t}|h(x)|.$$ Let us also consider the measures $M^p,\overline{M}^p$ associated to the perfect kernel $K_p$. 
Let us denote by $B_\lambda$ the ball centered at $0$ with radius $\lambda$. From Kahane's concentration inequalities \cite{Kah}, we have for all $q\leq 1$:
\begin{align*}
\E[(M(B_\lambda))^q] &\geq \E\big[\big(M^p(B_\lambda)e^{\gamma\sqrt{G_\lambda}Z-\frac{\gamma^2}{2 }G_\lambda}\big)^q\big]  
\end{align*}
where $Z$ is a standard Gaussian random variable independent of $M^p$. Hence, by using Theorem \ref{perfect}, we have:
\begin{align*}
\E[(M(B_\lambda))^q] &\geq \E\big[\big(M^p(B_\lambda)\big)^q\big]\E\big[\big(e^{\gamma\sqrt{G_\lambda}Z-\frac{\gamma^2}{2 }G_\lambda}\big)^q\big]  \\
&=\lambda^{\xi(q)}\E[M^p(  B_1 )^q]e^{q^2\frac{\gamma^2}{2 }G_\lambda-q\frac{\gamma^2}{2 }G_\lambda}.
\end{align*}
With the same argument we prove
\begin{align*}
e^{q^2\frac{\gamma^2}{2  }G_\lambda-q\frac{\gamma^2}{2  }G_\lambda}\E[(M(B_\lambda))^q] &\leq \lambda^{\xi(q)}\E[M^p(B_1 )^q].
\end{align*}
Because $G_\lambda\to 0$ as $\lambda\to 0$, the result follows from relation \eqref{relmom}.\qed

\section{Proofs of Section \ref{sec:KPZ}.} 
\subsection{Proof of Proposition \ref{prop:HSD} and Corollary \ref{cor:HSD}}
\noindent {\it Proof of Proposition \ref{prop:HSD}.} We assume that $A$ is bounded, say included in the ball $B(0,1)$.
We have for $ s<t$: 
$$H^{t,\delta}_\mu(A)\leq H^{s,\delta}_\mu(A)\sup_{\substack{B \: \text{ball}, \: B \cap A \not = \emptyset,\\B\subset B(0,1),{\rm diam}(B)\leq \delta}}\mu(B)^{t-s}.$$
Obviously, it suffices to prove that the above supremum converges to $0$ as $\delta\to 0$. This convergence is clear if  $\mu$ is diffuse. 

Let us now investigate the situation when  $\mu$ possesses atoms.  Let $A$ be a compact subset included in the ball $B(0,1)$, which does not encounter the atoms of $\mu$. We will prove:
$$\sup_{\substack{B \: \text{ball}, \: B \cap A \not = \emptyset,\\B\subset B(0,1),{\rm diam}(B)\leq \delta}}\mu(B) \to 0 \text{ as }\delta\to 0.$$
We argue by contradiction. Assume that this quantity does not converge towards $0$. We can find $\epsilon>0$, a sequence $(x_n)_n$ of points in $A$ and a sequence of balls $(B(y_n,r_n)_n)$ of radius $r_n$ going to $0$ such that $|y_n-x_n| \leq r_n$ and $\mu(B(y_n,r_n))\geq \epsilon$. Even if it means extracting a subsequence, we may assume that the sequence $(x_n)_n$ converges towards $x\in A$. We deduce $\mu(\{x\})\geq \epsilon$. This means that $\mu$ possesses an atom on $A$. Contradiction.\qed

\vspace{3mm}
\noindent {\it Proof of Corollary \ref{cor:HSD}.} If  $\mu=M$, then $\mu$ is diffuse (see Lemma \ref{atom} below). 
It remains to investigate the situation when $\mu=\overline{M}$. Let $A$ be a compact subset included in the ball $B(0,1)$ with null Lebesgue measure. For $0<\beta<\alpha$, we have
$\E[\overline{M}(A)^\beta]=c_{\alpha,\beta}\E[M(A)^\frac{\beta}{\alpha}]=0$ since $M(A)=0$ almost surely. Therefore, almost surely, the set $A$ does not encounter the atoms of $\overline{M}$.\qed

\begin{lemma}\label{atom}
Almost surely, the measure $M$ does not possess any atom.
\end{lemma}

\noindent {\it Proof.} By stationarity, it is enough to prove that, almost surely, the measure $M$ does not possess any atom on the cube  $[0,1]^d$. For $n\in\N^*$ and $k_1,\dots,k_d\in \{1,\dots,n\}$, let us denote by $I^n_{k_1,\dots,k_d}$ the cube $\prod_{i=1}^d[\frac{k_i-1}{n},\frac{k_i}{n}]$. From \cite[Corollary 9.3 VI]{daley}, it is enough to check that for each $\eta>0$:
$$\sum_{k_1,\dots,k_d=1}^n\P\Big(M(I^n_{k_1,\dots,k_d})>\eta\Big) = n^d\P\Big(M(I_{0,\dots,0}^n)>\eta\Big)\to 0\quad \text{ as }n\to \infty.$$
This is a direct consequence of the Markov inequality 
$$n^d\P\Big(M(I_{0,\dots,0}^n)>\eta\Big)\leq \frac{n^d}{\eta^{q} }\E[M(I_{0,\dots,0}^n)^{q}]$$
and the relation, for $1<q<\frac{2d}{\gamma^2}$ (see the proof of corollary \ref{notperf}),
$$\E[M(I_{0,\dots,0}^n)^{q}]\leq C n^{-\xi(q)}.$$
Indeed, for $1<q<\frac{2d}{\gamma^2}$, we have $\xi(q)>d$.\qed

\subsection{Proof of the standard KPZ formula}
The usual KPZ relation has already been proved in \cite{cf:DuSh} in the case of the $2d$-Gaussian Free Field in terms of expected box counting dimensions and in \cite{cf:RhoVar} for log-infinitely multifractal random measures in any dimensions in terms of almost sure Hausdorff dimensions. Let us stress that 
log-infinitely multifractal random measures strictly contains the class of Gaussian multiplicative chaos, which itself strictly contains the case of the GFF. For the sake of clarity and completeness, we sketch here a simple proof in the Gaussian case, which is based on the paper \cite{Benj}. The main differences  are the following:
\begin{enumerate}
\item We will use here exact scaling relations introduced in \cite{bacry}, and generalized to higher dimensions in \cite{rhovar}.
\item The correlations of the field are stronger than in \cite{Benj}, where independence is prominent.
\item We give a definition of Hausdoff dimension in terms of underlying measure which can be used in dimensions higher than one (the metric based definition of Hausdorff dimension used  \cite{Benj} is stuck to the dimension $1$ as long as one cannot give a proper definition of the Liouville metric).
\end{enumerate}
 
We further have the feeling that this short proof is worth being written as it helps to understand the KPZ formula in an easy way. It relies on the intensive use of the scaling properties of the Gaussian multiplicative chaos as well as the use of changes of probability measures (of Girsanov's type), which much simplifies the computations in comparison with \cite{cf:DuSh,cf:RhoVar}. For the sake of simplicity of notations, we make the proof in dimension $d=1$ but the proof in higher dimensions can be identically reproduced with minor modifications. We also assume that $M $ is the perfect measure, namely the measure with associated kernel given by $\gamma^2\ln_+\frac{T}{|x|}$.  Actually, it can easily be proved with the Kahane convexity inequalities (see \cite{Kah} or \cite[cor. 6.2]{cf:RoVa}) that this is not a restriction. We also mention that $M$ can be constructed as the limit
$$M(dx)=\lim_{l\to 0}M_l(dx)\stackrel{def}{=}e^{\gamma X^l_x-\frac{\gamma^2}{2}\E[(X_x^l)^2]}\,dx$$ where $X_l$ is a stationary Gaussian process with covariance kernel given by:
$$k_l(x)=\left\{\begin{array}{ll}
0 &\text{if }|x|>T,\\
 \ln_+\frac{T}{|x|} & \text{if }l T\leq |x|\leq T,\\
 \ln \frac{1}{l}+ \big(1- \frac{|x|}{Tl}\big)                &\text{if } 0\leq |x|\leq lT.
 \end{array}\right.$$
 Observe that such kernels fit the formalism of Section \ref{back:kah} by taking $l=2^{-n}$ for instance since we have
 $$M(dx)=\lim_{l\to 0}M_l(dx)=\lim_{n\to 0}M_{2^{-n}}(dx).$$
 Such a family of kernels possesses useful scaling properties, namely that for $|x|\leq T$ and $0<\lambda <1$, $k_{\lambda l}(\lambda x)=k_l(x)+\ln \frac{1}{\lambda}$. In particular, we have the following scaling relation for  all $0<l<1$ and all $0<\lambda<1$:
\begin{equation}\label{scalingkl}
 \big((X^{\lambda l}_{\lambda x})_{x\in B(0,T)},(M_{\lambda l}(\lambda A))_{A\subset B(0,T)}\big)\stackrel{law}{=}\big((X^l_{ x}+\Omega_\lambda)_{x\in B(0,T)},(\lambda e^{\gamma\Omega_\lambda-\frac{\gamma^2}{2}\ln\frac{1}{\lambda}}M_{  l}(  A))_{A\subset B(0,T)}\big).
\end{equation}
where $\Omega_\lambda$ is a centered Gaussian random variable with variance $\ln\frac{1}{\lambda}$ and independent of the couple $\big((X^l_{ x} )_{x\in B(0,T)},( M_{  l}(  A))_{A\subset B(0,T)}\big)$. We will use the above relation throughout the proof.

\vspace{2mm}
Now we begin with the proof. Without loss of generality we assume that $T=1$. Let $K$ be a compact subset of $\R $, included in $[0,1]$, with Hausdorff dimension $ {\rm dim}_{Leb}(K)$. Let $q\in [0,1]$ be such that $\xi(q)>{\rm dim}_{Leb}(K)$.
For $\epsilon>0$, there is a covering of $K$ by a countable family of balls $(B(x_n,r_n))_n$ such that 
$$ \sum_n r_n^{\xi(q)}<\epsilon.$$ By using in turn  the stationarity and the power law spectrum of the measure, we have
\begin{align*}
\E\Big[\sum_nM(B(x_n,r_n))^q\Big]& =\sum_n\E\Big[M(B(0,r_n))^q\Big]\\
&\leq C_q \sum_n r_n^{\xi(q)}\\
&\leq C_q\epsilon,
\end{align*}
we deduce by the Markov inequality
$$\P\Big(\sum_nM(B(x_n,r_n))^q\leq C_q \sqrt{\epsilon}\Big)\geq 1-\sqrt{\epsilon}.$$
Thus, with probability $1-\sqrt{\epsilon}$, there is a covering of balls of $K$ such that 
$\sum_nM(B(x_n,r_n))^q\leq C_q\sqrt{\epsilon}$. So $q\geq {\rm dim}_M(K)$ almost surely.

Conversely, consider $q\geq 0$ such that $\xi(q)<{\rm dim}_{Leb}(K)$. By the Frostman Lemma, there is a probability measure $\kappa$ supported by $K$ such that
$$\int_{[0,1]^2}\frac{1}{|x-y|^{\xi(q)}}\kappa(dx)\kappa(dy)<+\infty.$$ Let us define the random measure $\widetilde{\kappa}$ as the almost sure limit of the following family of positive random measures:
\begin{equation}\label{gamtilde}
\widetilde{\kappa}(dx)=\lim_{l\to 0}e^{q\gamma X^l_x-\frac{q^2\gamma^2}{2}\E[(X^l_x)^2]}\kappa(dx).
\end{equation} 
For $\gamma^2<2$, the limit is non trivial because $q^2\gamma^2/2<\xi(q)$ and supported by $K$  (see \cite{Kah}). From the Frostman lemma again, we just have to prove that the quantity
\begin{equation}\label{espfrost}
\int_{[0,1]^2}\frac{1}{M([x,y])^{q}}\widetilde{\kappa}(dx)\widetilde{\kappa}(dy)
\end{equation} 
is finite almost surely. Actually, since we use a non-standard version of the Frostman lemma (in terms of measures, not distances), we provide a proof in the next lemma. 

Coming back to the proof of KPZ, it suffices to prove that the quantity \eqref{espfrost} has a finite expectation. Moreover, by using the Fatou lemma and the stationarity of the measure $M$,  we have
\begin{align*}
\E\Big[\int_{[0,1]^2}\frac{1}{M([x,y])^{q}}\widetilde{\kappa}(dx)\widetilde{\kappa}(dy)\Big]&\leq \liminf_l \int_{[0,1]^2}\E\Big[\frac{e^{q\gamma X^l_x+q\gamma X^l_y- q^2\gamma^2\E[(X^l_x)^2]}}{M_l([x,y])^{q}}\Big] \kappa(dx) \kappa (dy)\\
&= \liminf_l 2\int_{y\geq x}\E\Big[\frac{e^{q\gamma X^l_0+q\gamma X^l_{y-x} -q^2\gamma^2\E[(X^l_x)^2]}}{M_l([0,y-x])^{q}}\Big] \kappa(dx) \kappa (dy).
 \end{align*}
 We decompose the last integral into two terms:
\begin{align*}
\int_{y\geq x}\E\Big[&\frac{e^{q\gamma X^l_0+q\gamma X^l_{y-x} -q^2\gamma^2\E[(X^l_x)^2]}}{M_l([0,y-x])^{q}}\Big] \kappa(dx) \kappa (dy)\\=&
\int_{0\leq y-x\leq l}\E\Big[\frac{e^{q\gamma X^l_0+q\gamma X^l_{y-x} -q^2\gamma^2\E[(X^l_x)^2]}}{M_l([0,y-x])^{q}}\Big] \kappa(dx) \kappa (dy)\\&+\int_{y- x\geq l}\E\Big[\frac{e^{q\gamma X^l_0+q\gamma X^l_{y-x} -q^2\gamma^2\E[(X^l_x)^2]}}{M_l([0,y-x])^{q}}\Big] \kappa(dx) \kappa (dy)\\
 \stackrel{def}{=}&A^1_l+A_l^2.
\end{align*}
For each of the above terms, we will use an appropriate scaling relation.

By using \eqref{scalingkl}, we deduce
\begin{align*}
A^2_l&=    \int_{y- x\geq l}\E\Big[\frac{e^{2q\gamma \Omega_{y-x}-q^2\gamma^2\ln\frac{1}{y-x}}e^{q\gamma X^{\frac{l}{y-x}}_0+q\gamma X^{\frac{l}{y-x}}_{1} -q^2\gamma^2\E[(X^{\frac{l}{y-x}}_x)^2]}}{(y-x)^qe^{q\gamma\Omega_{y-x}-q\frac{\gamma^2}{2}\ln\frac{1}{y-x}}M_{\frac{l}{y-x}}([0,1])^q}\Big] \kappa(dx) \kappa (dy)\\
&= \int_{y- x\geq l}\E\Big[\frac{e^{ q\gamma \Omega_{y-x}-(q^2\gamma^2-q\frac{\gamma^2}{2})\ln\frac{1}{y-x}}}{(y-x)^q}\Big]\E\Big[\frac{e^{q\gamma X^{\frac{l}{y-x}}_0+q\gamma X^{\frac{l}{y-x}}_{1} -q^2\gamma^2\E[(X^{\frac{l}{y-x}}_x)^2]}}{ M_{\frac{l}{y-x}}([0,1])^{q}}\Big] \kappa(dx) \kappa (dy)\\
&=   \int_{y- x\geq l}\frac{1}{(y-x)^{\xi(q)}} \E\Big[\frac{e^{q\gamma X^{\frac{l}{y-x}}_0+q\gamma X^{\frac{l}{y-x}}_{1} -q^2\gamma^2\E[(X^{\frac{l}{y-x}}_x)^2]}}{ M_{\frac{l}{y-x}}([0,1])^{q}}\Big] \kappa(dx) \kappa (dy)
 \end{align*}
By using a change of measures, we have
\begin{align*}
\E\Big[\frac{e^{q\gamma X^{\frac{l}{y-x}}_0+q\gamma X^{\frac{l}{y-x}}_{1} -q^2\gamma^2\E[(X^{\frac{l}{y-x}}_x)^2]}}{ M_{\frac{l}{y-x}}([0,1])^{q}}\Big] &= \E\Big[\frac{e^{q^2\gamma^2k_{\frac{l}{y-x}}(1)}}{ \Big(\int_0^1e^{\gamma X^{\frac{l}{y-x}}_r-\frac{\gamma^2}{2}\E[(X^{\frac{l}{y-x}}_x)^2]+q\gamma^2k_{\frac{l}{y-x}}(1-r)+q\gamma^2k_{\frac{l}{y-x}}(r)}\,dr\Big)^q}\Big]  \\
&\leq C\E\Big[\frac{1}{ \Big(\int_{0}^{1}e^{\gamma X^{\frac{l}{y-x}}_r-\frac{\gamma^2}{2}\E[(X^{\frac{l}{y-x}}_x)^2]}\,dr\Big)^q}\Big]
\end{align*}
for some positive constant $C$. Notice that we have just used the fact that $k_{\frac{l}{y-x}}(1)=0$ and that $k_{\frac{l}{y-x}}$ is positive. 
It is a standard fact that the measure $M$ possesses moments of negative order (see \cite{Bar}) so that we have proved
$$\lim_l A^2_l\leq C\int_{[0,1]^2}\frac{1}{|y-x|^{\xi(q)}}\kappa(dx)\kappa(dy)<+\infty.$$
To treat the term $A^1_l$, we use quite a similar argument excepted that we use the scaling relation on $l$ instead of $y-x$, and a change of measures again:
\begin{align*}
A^1_l&=    \int_{0\leq y-x\leq l}\E\Big[\frac{e^{2q\gamma \Omega_{l}-q^2\gamma^2\ln\frac{1}{l}}e^{q\gamma X^{1}_0+q\gamma X^{1}_{\frac{y-x}{l}} -q^2\gamma^2\E[(X_x^1)^2]}}{l^qe^{q\Omega_{l}-q\frac{\gamma^2}{2}\ln\frac{1}{l}}M_{1}([0,\frac{y-x}{l}])^{q}}\Big] \kappa(dx) \kappa (dy)\\
&= \int_{0\leq y-x\leq l}\E\Big[\frac{e^{ q\gamma \Omega_{l}-(q^2\gamma^2-q\frac{\gamma^2}{2})\ln\frac{1}{l}}}{l^q}\Big]\E\Big[\frac{e^{q\gamma X_{\frac{y-x}{l}}^1+q\gamma X_0^{1} -q^2\gamma^2\E[(X^1_x)^2]}}{ M_{1}([0,\frac{y-x}{l}])^{q}}\Big] \kappa(dx) \kappa (dy)\\
&=   \int_{0\leq y-x\leq l}\frac{1}{l^{\xi(q)}} \E\Big[\frac{e^{q^2\gamma^2k_1( \frac{y-x}{l}})}{ \Big(\int_0^{\frac{y-x}{l}}e^{\gamma X_{r}^1-\frac{\gamma^2}{2}\E[(X^1_r)^2]+q\gamma^2k_1(\frac{y-x}{l}-r)+q\gamma^2k_1(r)}\,dr\Big)^q}\Big] \kappa(dx) \kappa (dy).
 \end{align*}
By using the fact that $k_1$ is positive and bounded by $1$, we have (for some positive constant $C$ independent of $l$)
$$A^1_l\leq C\int_{0\leq y-x\leq l}\frac{1}{l^{\xi(q)}} \E\Big[\frac{1}{ \Big(\int_0^{\frac{y-x}{l}}e^{\gamma X_{r}^1-\frac{ \gamma^2}{2}\E[(X^1_r)^2]}\,dr\Big)^q}\Big] \kappa(dx) \kappa (dy).$$
Since $\E[X^1_rX^1_0]\leq \E[(X^1_0)^2]$, we can use Kahane's convexity inequalities to the convex mapping $x\mapsto \frac{1}{x^q}$. We deduce (for some positive constant $C'$)
\begin{align*}A^1_l\leq &C\int_{0\leq y-x\leq l}\frac{1}{l^{\xi(q)}} \E\Big[\frac{1}{ \Big(\int_0^{\frac{y-x}{l}}e^{\gamma X_{0}^1-\frac{ \gamma^2}{2}\E[(X^1_0)^2]}\,dr\Big)^q}\Big] \kappa(dx) \kappa (dy)\\
\leq & C'\int_{0\leq y-x\leq l}\frac{l^q}{l^{\xi(q)}(y-x)^q}  \kappa(dx) \kappa (dy)\\
\leq &C'\int_{0\leq y-x\leq l}\frac{1}{(y-x)^{\xi(q)} }  \kappa(dx) \kappa (dy).
\end{align*}
Hence $$\lim_l A^1_l\leq C'\int_{B(0,T)^2}\frac{1}{|y-x|^{\xi(q)}}\kappa(dx)\kappa(dy)<+\infty.$$
The KPZ formula is proved (by using scaling relations only).\qed

Now we give a proof of the Frostman lemma. We deal with dimension greater than $1$ so that the reader has no difficulty to generalize the proof of the KPZ formula to dimensions higher than $1$. For any $x\in\R^d$ and $r>0$, the ball $B(x,r)$ can be decomposed into exactly $2^d$ isometric pieces that can be written as
$$B(x,r)\cap\Big\{z=(z_1,\dots,z_d)\in\R^d; (z_1-x_1) \epsilon_1\geq 0,\dots, (z_d-x_d) \epsilon_d\geq 0\Big\} $$ where $\epsilon_1,\dots,\epsilon_d\in\{-1,1\}$ (the fact that they may overlap at their boundary does not matter). Let us enumerate by $(HB(x,r,i))_{ i \in\{-1;1\}^d}$ these portions of balls. 

\begin{lemma}[Frostman lemma]
Assume that $\mu,\nu$ are two Radon measures on $\R^d$. Assume further that $\nu$ is a probability measure supported by a compact set $K\subset B(0,1)$. If for all $i\in \{-1,1\}^d$,
$$\int_K\int_K\frac{\nu(dx)\nu(dy)}{ \mu\Big(HB(\frac{x+y}{2},\frac{|x-y|}{2},i)\Big)^q}<+\infty$$ then ${\rm dim}_\mu(K)\geq q$.
\end{lemma}

\noindent {\it Proof.} Let us define the function:
$$\forall x\in K,\quad g(x)=\int_K\frac{\nu(dy)}{\inf_{i\in\{-1,1\}^d}\mu\Big(HB(\frac{x+y}{2},\frac{|x-y|}{2},i)\Big)^q}.$$
Observe that the assumptions imply that $\int_Kg(x)\nu(dx)< \infty$. We deduce that 
$$\nu\big(\{x\in K;g(x)\leq L\}\big)\to 1,\quad \text{ as }L\to\infty.$$ Therefore we can find $L$ large enough such that the set $K_L=\{x\in K;g(x)\leq L\}$ satisfies $\nu(K_L)\geq \frac{1}{2}$. Let us consider a covering $(B(x_n,r_n))_n$ of $K$ with balls of radius less than $\delta$. We consider the subsequence $(B(x_{n_k},r_{n_k}))_{n_k}$ of balls which intersect $K_L$. It is obvious that this subsequence is a covering of $K_L$. For each $n_k$, there exists $y_{n_k}$ in $K_L \cap B(x_{n_k},r_{n_k})$. From Lemma \ref{lem.incl} below, for all $y \in B(x_{n_k},r_{n_k})$, one of the portion of balls $HB\big(\frac{y_{n_k}+y}{2},\frac{|y_{n_k}-y|}{2},i\big)$ for $i\in \{-1,1\}^d$ is included in $B(x_{n_k},r_{n_k})$. Hence: 
\begin{equation*}
\inf_{i\in\{-1,1\}^d} \mu\Big(HB\big(\frac{y_{n_k}+y}{2},\frac{|y_{n_k}-y|}{2},i\big)\Big)^q  \leq \mu( B(x_{n_k},r_{n_k})  )^q.
\end{equation*}
Therefore, we get:
\begin{equation*}
\frac{\nu(  B(x_{n_k},r_{n_k}) )}{\mu(  B(x_{n_k},r_{n_k}) )^q  } \leq  \int_K\frac{\nu(dy)}{\inf_{i\in\{-1,1\}^d} \mu\Big(HB\big(\frac{y_{n_k}+y}{2},\frac{|y_{n_k}-y|}{2},i\big)\Big)^q}   \leq L
\end{equation*} 
This leads to 
\begin{equation*}
\sum_k  \mu(  B(x_{n_k},r_{n_k}) )^q  \geq \frac{1}{L} \sum_k  \nu(  B(x_{n_k},r_{n_k}) ) \geq \frac{1}{L} \nu(  K_L)  \geq \frac{1}{2L}
\end{equation*}
which gives the result.\qed

\begin{lemma}\label{lem.incl}
Consider $x,y,z\in\R^d$ and $r>0$ such that $x,y\in B(z,r)$. Then, for some $i\in\{-1,1\}^d$, 
$$HB\big(\frac{x+y}{2},\frac{|y-x|}{2},i\big)\subset B(z,r).$$
\end{lemma}

\noindent {\it Proof.} If $x=y$, the proof is obvious. So we assume $x\not =y$.  Let us set $c=(x+y)/2$ and  assume that $c\not =z$, otherwise the proof is trivial. Let us then define the half-space
$$H=\{u\in\R^d;(u-c,z-c)\geq 0\}.$$
It is plain to check that the half-ball $H\cap B(c,\frac{|x-y|}{2})$ is contained in $B(z,r)$. Then, it is also elementary to notice that any non trivial half-space whose boundary contains some point  $v\in\R^d$ necessarily contains at least one of the portion of space $$\Big\{z=(z_1,\dots,z_d)\in\R^d; (z_1-v_1) \epsilon_1\geq 0,\dots, (z_d-v_d) \epsilon_d\geq 0\Big\} $$ where $\epsilon_1,\dots,\epsilon_d\in\{-1,1\}$. 
Indeed, if this were not true, the complementary of such a half-space should contain an interior point of each of these portions, and therefore should contain an open neighborhood of $v$ by convexity, contradiction.
Therefore, we can find $\epsilon_1,\dots,\epsilon_d\in\{-1,1\}$ such that
$$\Big\{z=(z_1,\dots,z_d)\in\R^d; (z_1-c_1) \epsilon_1\geq 0,\dots, (z_d-c_d) \epsilon_d\geq 0\Big\} \cap B(c,\frac{|x-y|}{2})\subset B(z,r),$$
which complete the proof.\qed

\subsection{Proof of the dual KPZ formula} \label{proof:dual}
This time, we do not restrict to the dimension $1$. Let $K$ be a compact subset of $\R^d $, included in the ball $B(0,1)$ with Hausdorff dimension $0\leq {\rm dim}_{Leb}(K)<1$.  Let $\delta_0 $ be the unique solution in $[0,\alpha[$ such that $\frac{\overline{\xi}(\delta_0)}{d}={\rm dim}_{Leb}(K)$. We want to prove $\delta_0={\rm dim}_{\overline{M}}(K)$.

Let $0\leq q<\alpha $ be such that $\frac{\overline{\xi}(q)}{d}>{\rm dim}_{Leb}(K)$.
For $\epsilon>0$, there is a covering of $K$ by a countable family of balls $(B(x_n,r_n))_n$ such that 
$$ \sum_n r_n^{ \overline{\xi}(q)}<\epsilon.$$ Since we have (see Theorem \ref{perfect})
\begin{align*}
\E\Big[\sum_n\overline{M}(B(x_n,r_n))^q\Big]& =\sum_n\E\Big[\overline{M}(B(0,r_n))^q\Big]\\
&\leq C_q \sum_n r_n^{\overline{\xi}(q)}\\
&\leq C_q\epsilon,
\end{align*}
we deduce by the Markov inequality
$$\P\Big(\sum_n\overline{M}(B(x_n,r_n))^q\leq C_q \sqrt{\epsilon}\Big) \geq 1-\sqrt{\epsilon}.$$
Thus, with probability $1-\sqrt{\epsilon}$, there is a covering of balls of $K$ such that 
$\sum_n\overline{M}(B(x_n,r_n))^q\leq C_q\sqrt{\epsilon}$. So $q\geq {\rm dim}_{\overline{M}}(K)$ almost surely.

Conversely, consider $p\in [0,\alpha[$ such that $\frac{\overline{\xi}(p)}{d}<{\rm dim}_{Leb}(K)$. Since $\overline{\xi}(p)=\xi(\frac{p}{\alpha})$, we can set $q=\frac{p}{\alpha}\in  [0,1[$ and we have $\frac{\xi(q)}{d}<{\rm dim}_{Leb}(K)$. As we proved above, we can consider the measure $\widetilde{\kappa}$ introduced in \eqref{gamtilde}. It is almost surely supported by $K$ and non trivial. Furthermore, it satisfies  
$$\E\Big[\int_{B(0,T)^2}\frac{1}{M(B(x,|y-x|))^{q}}\widetilde{\kappa}(dx)\widetilde{\kappa}(dy)\Big]<+\infty.$$
Let us prove that
\begin{equation}\label{newproof}
\E\Big[\int_{B(0,T)^2}\frac{1}{\overline{M}(B(x,|y-x|))^{p}}\widetilde{\kappa}(dx)\widetilde{\kappa}(dy)\Big]<+\infty.
\end{equation}
By using the relation for $p,x>0$
$$\Gamma(p)=x^p\int_0^{+\infty}u^{p-1}e^{-ux}\,du,$$
we deduce:
\begin{align*}
\E\Big[\int_{B(0,T)^2}&\frac{1}{\overline{M}(B(x,|y-x|))^{p}}\widetilde{\kappa}(dx)\widetilde{\kappa}(dy)\Big] \\
&=\frac{1}{\Gamma(p)}\E\Big[\int_0^{+\infty}u^{p-1}\int_{B(0,T)^2}e^{-u \overline{M}(B(x,|y-x|))}\widetilde{\kappa}(dx)\widetilde{\kappa}(dy)\,du\Big]\\
&=\frac{1}{\Gamma(p)}\E\Big[\int_0^{+\infty}u^{p-1}\int_{B(0,T)^2}\E\Big[e^{-u \overline{M}(B(x,|y-x|))}|Y^n,n\geq 1\Big]\widetilde{\kappa}(dx)\widetilde{\kappa}(dy)\,du\Big]\\
&=\frac{1}{\Gamma(p)}\E\Big[\int_0^{+\infty}u^{p-1}\int_{B(0,T)^2} e^{-u^\alpha M(B(x,|y-x|))}\widetilde{\kappa}(dx)\widetilde{\kappa}(dy)\,du\Big]
\end{align*}
Now we make the change of variables $y=u^\alpha M(B(x,|y-x|))$ to obtain:
\begin{align*}
\E\Big[\int_{B(0,T)^2}&\frac{1}{\overline{M}(B(x,|y-x|))^{p}} \widetilde{\kappa}(dx)\widetilde{\kappa}(dy)\Big] \\
&= \frac{1}{\alpha\Gamma(p)}\E\Big[\int_{B(0,T)^2}\frac{1}{M(B(x,|y-x|))^q}  \widetilde{\kappa}(dx)\widetilde{\kappa}(dy)\Big]\int_0^{+\infty}y^{\frac{p}{\alpha}-1}e^{-y}\,dy\\
&=\frac{\Gamma(\frac{p}{\alpha}+1)}{ \Gamma(p+1)}\E\Big[\int_{B(0,T)^2}\frac{1}{M(B(x,|y-x|))^q}  \widetilde{\kappa}(dx)\widetilde{\kappa}(dy)\Big] .
\end{align*}
Hence, the above quantity is finite and \eqref{newproof} is proved. In fact, with minor modifications, one can prove: 
\begin{equation*}
\E\Big[\int_K\int_K\frac{\widetilde{\kappa}(dx)\widetilde{\kappa}(dy)}{\overline{M}\Big(HB(\frac{x+y}{2},\frac{|x-y|}{2},i)\Big)^q}\Big]<+\infty.
\end{equation*}
for $i\in\{-1,1\}^d$. 
We conclude by using the above Frostman lemma that $p<{\rm dim}_{\overline{M}}(K)$. The dual KPZ formula is proved. Notice that we have also proved  the relation $\dim_{\overline{M}}(K)=\alpha \dim_M(K)$, which is nothing but the duality relation. Finally, we stress that the argument for the dual KPZ formula can be obviously generalized and works for any measure $M$ and its subordinated counterpart.  \qed
 
\section{Proof of Remark \ref{rem.otherlaw}}\label{otherlaw}
Let $\mathcal{X}=(X_j)$ be a collection of compact subsets of $\R^d$ with disjoint interiors such that $\bigcup_j X_j=\R^d$, and $(Z_k)$ a partition of $\mathbb R_+^*$ into compact subintervals semi-open to the right. Let $\mathcal  M=\mathcal{M}(\mathcal{X})$ be the collection of Radon measures on $\R^d$ satisfying the following property: 
\[
  \mu(\partial X_j)=0 \text{ for each } j.
\]
Endow $\mathcal{M}$ with the topology of weak convergence.  Let $\mathcal  M'$ be the space of probability measures on 
$$
Y=\Big (\prod_{j,k} (\mathbb N\times (X_j\times Z_k)^{\mathbb{N}_+}), \otimes_{j,k} (\mathcal B(\mathbb{N})\otimes (\mathcal B(X_j)\otimes \mathcal B(Z_k))^{\otimes \mathbb{N_+}}) \Big).
$$
This space can be endowed with the Wasserstein distance of order 1 \cite[Ch. 6]{Villani} ($Y$ being endowed with its natural structure of Polish space). 

Denote by $\rho_\alpha$ the measure $dz/z^{1+\alpha}$ over $\mathbb R_+^*$. Denote by $\mathcal  P(\lambda)$ the Poisson distribution with parameter $\lambda$. If $M\in\mathcal  M$, let $\nu(M)$ be the element of $\mathcal  M'$ defined as 
$$\otimes_{j,k} \mathcal P(M(X_j)\rho_\alpha(Z_k))\otimes \left(\frac{M_{|X_j}}{M(X_j)}\otimes \frac{{\rho_\alpha}_{|Z_k}}{\rho_\alpha(Z_k)} \right)^{\otimes \mathbb{N}_+}.$$
The mapping $M\mapsto \nu(M)$ is continuous. Then consider the continuous mapping $N: Y\to \mathcal M$ defined as 
$$
N(y=(n_{j,k},((x_{j,i},z_{k,i}))_{i\ge 1})_{j,k})=\sum_{j,k}\sum_{i=1}^{n_{j,k}}z_{j,i}\delta_{(x_{j,i},z_{k,i})}.
$$

We work with the limit multiplicative chaos $M=M(\omega)$ with the $d$-dimensional Lebesgue measure $\mathcal{L}$ as the reference measure. One may choose a partition $\mathcal{X}$ such that $\mathcal{L}\in\mathcal{M}(\mathcal{X})$, for example
\[
\mathcal{X}=\left\{\prod_{k=1}^d[n_k,n_k+1]: n_k\in \mathbb{Z}\right\}.
\]
Now using the countability of $\mathcal X$ and the fact that for each fixed Borel subset $A$ of $\R^d$, $\mathcal L(A)=0$ implies $M(A)= 0$ almost surely, we see that  $M\in \mathcal{M}(\mathcal{X})$ almost surely. Then the random measure $(\omega,y)\mapsto N(y)$ defined on $\Omega\times Y$ endowed with the skew product measure $ \nu(M(\omega))(dy)\mathbb P(d\omega)$ provides a measurable construction for the law of the Poisson random measure as suggested in Remark \ref{rem.otherlaw}. 
\hspace{10 cm}

\end{document}